\documentclass[11pt]{article}

\usepackage[T1]{fontenc}
\usepackage[a4paper,margin=29mm]{geometry}
\usepackage{amsmath,amssymb,amsthm,mathtools,mathrsfs}
\usepackage{enumitem}
\usepackage{microtype}
\usepackage{hyperref}
\usepackage{bm}
\usepackage[title]{appendix}

\hypersetup{colorlinks=true,linkcolor=blue,citecolor=blue,urlcolor=blue,
  pdftitle={Calculus of Variations for Degenerate Hamiltonian Actions},
  pdfauthor={Takaaki Nishida and Kohei Soga}}

\newtheorem{theorem}{Theorem}[section]
\newtheorem{proposition}[theorem]{Proposition}
\newtheorem{lemma}[theorem]{Lemma}
\newtheorem{corollary}[theorem]{Corollary}
\theoremstyle{definition}

\theoremstyle{plain}
\newtheorem{remark}[theorem]{Remark}

\newcommand{\T}{\mathbb T}
\newcommand{\R}{\mathbb R}
\newcommand{\eps}{\varepsilon}
\newcommand{\dd}{\,\mathrm d}
\newcommand{\Lip}{\operatorname{Lip}}

\newcommand{\Ac}{\mathcal A^c}
\newcommand{\AcK}{\mathcal A_K^c}

\numberwithin{equation}{section}

\title{Calculus of Variations for Degenerate Hamiltonian Actions}
\author{Takaaki Nishida\thanks{Graduate School of Informatics, Kyoto University \newline Yoshida-Honmachi, Sakyo-ku, Kyoto, 606-8501 Japan. \quad Email: tkknish@acs.i.kyoto-u.ac.jp} \and Kohei Soga\thanks{Department of Mathematics, Faculty of Science and Technology, Keio University \newline 3-14-1 Hiyoshi, Kohoku-ku, Yokohama, 223-8522 Japan. \quad Email: soga@math.keio.ac.jp}}
\date{}

\begin{document}
\maketitle

\begin{abstract}
We study one-dimensional variational problems for the Hamiltonian action generated by a Hamiltonian
$H\in C^2(\mathbb T\times\mathbb T\times[0,C];\mathbb R)$
that is convex, but not necessarily strictly convex, in the momentum variable $p\in[0,C]$.
We also assume that $H(x,t,0)$ and $H(x,t,C)$ are independent of $x$.
The standard Tonelli theory for the Lagrangian action obtained by the Legendre transform of $H$ need not apply, because the corresponding Lagrangian need not be differentiable.
We therefore work directly with the Hamiltonian action and formulate an  $\inf$-$\sup$ problem for optimal trajectories.
For the problem with two fixed endpoints, we prove that every reachable pair of endpoints admits a $C^2$-smooth optimal trajectory satisfying Hamilton's equations.
This result holds even when \(H\) is degenerate in \(p\), in which case maximizing momentum curves associated with a minimizing curve need not be unique or smooth.
We also prove the corresponding result for the Bolza problem.
The proof uses a finite-dimensional approximation and a small quadratic regularization in the momentum variable.
One-sided variations control the momentum boundary and give the discrete form of Hamilton's equations.
These equations give uniform estimates for  the discrete optimal  trajectories,  yielding the strong compactness required to pass to the limit. 
A further comparison argument proves that the limiting Hamiltonian trajectory is optimal for the original  $\inf$-$\sup$ problem.
We also give a concrete physical model leading to this class of Hamiltonians.
The results provide a good starting point for the dynamical aspect of weak KAM theory in this degenerate setting.
\end{abstract}

\medskip
\noindent\textbf{Keywords.} Hamiltonian action; degenerate convex Hamiltonian;  $\inf$-$\sup$ variational problem; finite-dimensional approximation; one-sided variation.

\medskip
\noindent\textbf{2020 Mathematics Subject Classification.} Primary 49J35; Secondary 37J51, 49J15.

\section{Introduction}\label{sec:introduction}

We consider a time-periodic one-dimensional Hamiltonian whose momentum variable belongs to the compact interval $[0,C]$, where $C>0$.  We assume:
\begin{description}
\item[\textup{(H1)}] $H\in C^2(\T\times\T\times[0,C];\R)$, \qquad $\T:=\R/(2\pi\mathbb Z)$;
\item[\textup{(H2)}] $\partial_p^2 H(x,t,p)\ge0$ \quad for every $(x,t,p)\in\T\times\T\times[0,C]$;
\item[\textup{(H3)}] $H(x,t,0)=b_0(t)$, \qquad $H(x,t,C)=b_C(t)$.
\end{description}
We identify such an $H$ with its $2\pi$-periodic lift $H:\R\times\R\times[0,C]\to\R$ in the $x$- and $t$-variables.  The variational problems below are considered for this lift.
We assume only convexity in $p$, not strict or uniform convexity, and impose no special condition on the degeneracy set $\{(x,t,p)\in\T\times\T\times[0,C]\,|\,\partial_p^2H(x,t,p)=0\}$.  Thus complete degeneracy in the momentum variable is allowed. 

The class of Hamiltonians considered here is motivated by reductions of autonomous Hamiltonian systems with two degrees of freedom, restricted to fixed energy surfaces, to nonautonomous Hamiltonian systems with one degree of freedom. 
The invariance of the boundary of the momentum domain under the original Hamiltonian flow leads to (H3), after a suitable symplectic change of variables. 
Section 2 illustrates this mechanism through a concrete model based on the two-particle $\beta$-FPU oscillator,  leading to a Hamiltonian satisfying (H1)--(H3).

The boundary condition (H3) ensures that the momentum boundaries $p=0$ and $p=C$ are invariant under the Hamiltonian flow. This property plays an essential role in our analysis.

\subsection{The variational problem for the Hamiltonian action}\label{subsec:intro-problem} 
Fix $c\in(0,C)$ and a terminal time $t>0$.  For a Lipschitz configuration-space curve $\gamma:[0,t]\to\R$ and a measurable momentum curve $p:[0,t]\to[0,C]$, define\footnote{Here and below, a dot denotes differentiation with respect to $s$.}
\begin{equation*}%\label{eq:intro-action}
 \mathcal A^c(\gamma(\cdot),p(\cdot))
 :=\int_0^t
 \{\dot\gamma(s)(p(s)-c)-H(\gamma(s),s,p(s))\}\,\mathrm ds.
\end{equation*}
By setting $p^c(\cdot):=p(\cdot)-c$ and $H^c(x,t,q):=H(x,t,c+q)$, the action $\Ac$ takes the standard form of a Hamiltonian action.\footnote{Indeed, $p^c:[0,t]\to[-c,C-c]$ and
%\[
 $\displaystyle \mathcal A^c(\gamma(\cdot),p(\cdot))
 =\int_0^t\{\dot\gamma(s)p^c(s)-H^c(\gamma(s),s,p^c(s))\}\,\mathrm ds. $
%\]
}
We call the curve $\gamma(\cdot)$ {\it admissible} if
\begin{equation}\label{eq:intro-admissible}
 \partial_p H(\gamma(s),s,0)\le \dot\gamma(s)
 \le \partial_p H(\gamma(s),s,C)
 \quad\text{for a.e. }s.
\end{equation}
We impose this velocity constraint a priori, before taking variations. 
Define
\begin{align*}
\Gamma_t
&:=\left\{\gamma\in\Lip([0,t];\R)\,|\,\mbox{ $\gamma$ satisfies \eqref{eq:intro-admissible}}\right\},\\
P_t
&:=\left\{p\in L^\infty(0,t)\,|\,0\le p(s)\le C
 \text{ for a.e. }s\right\}.\notag
\end{align*}
A pair $(\gamma(\cdot),p(\cdot))\in\Gamma_t\times P_t$ is called an admissible trajectory.  
For two fixed endpoints $\gamma(0)=y$ and $\gamma(t)=x$, we consider
\begin{equation}\label{eq:intro-infsup}
 \inf_{\substack{\gamma(\cdot)\in \Gamma_t \\
                   \gamma(0)=y,\ \gamma(t)=x}}
 \ \sup_{p(\cdot)\in P_t}\mathcal A^c(\gamma(\cdot),p(\cdot)).
\end{equation}
The set of admissible curves over which the infimum in  \eqref{eq:intro-infsup} is taken is nonempty only for reachable initial points $y$ (note that we first fix an arbitrary $x$ and then choose $y$ so that the range is nonempty); the reachable set is characterized in Section~\ref{sec:reachable}.
For a fixed admissible curve $\gamma(\cdot)$, a momentum curve attaining the inner supremum in \eqref{eq:intro-infsup} is called a \emph{maximizing momentum curve associated with $\gamma(\cdot)$};  a curve attaining the outer infimum is called a \emph{minimizing curve}.  A pair $(\gamma^*(\cdot),p^*(\cdot))$ consisting of a minimizing curve and an associated maximizing momentum curve is called an \emph{optimal trajectory} for the Hamiltonian action $\Ac$.  This is only a variational definition, i.e., an optimal trajectory is not a priori assumed to satisfy Hamilton's equations.  We call \eqref{eq:intro-infsup} the \emph{two-fixed-endpoint  $\inf$-$\sup$ problem}.

The main problem of this paper is to seek an optimal trajectory that is also a classical Hamiltonian trajectory.  If $p(\cdot)$ is a maximizing momentum curve associated with an admissible curve $\gamma(\cdot)$, so that $p(\cdot)$  is a pointwise maximizer for almost every $s$\footnote{This means $\displaystyle\dot{\gamma}(s)(p(s)-c)-H(\gamma(s),s,p(s))=\max_{0\le p\le C}\{ \dot{\gamma}(s)(p-c)-H(\gamma(s),s,p) \}$ for almost every $s$.}, then the admissibility condition \eqref{eq:intro-admissible} and convexity of $H$ in $p$ give $\dot\gamma(s)=\partial_pH(\gamma(s),s,p(s))$ for almost every $s$.  This is the first part of Hamilton's equations.  The difficulty is to choose $p(\cdot)$ so that the second part of Hamilton's equations also holds and the resulting trajectory is $C^2$-smooth in the degenerate situation. 

We also consider the corresponding Bolza problem, that is, the problem with one fixed endpoint and an initial cost.  This is the variational problem naturally associated with the initial-value problem for the corresponding Hamilton--Jacobi equation.

The parameter $c$ is a reference momentum rather than a dynamical variable. From the viewpoint of weak KAM theory, $c$ is regarded as a cohomology parameter. In the two-fixed-endpoint problem the term $-c\dot\gamma(\cdot)$ contributes only $-c(x-y)$, while it affects the minimization in the Bolza problem because the initial endpoint is free. 

We remark that the  $\inf$-$\sup$ problem \eqref{eq:intro-infsup} also has an equivalent Lagrangian representation even when $H$ is degenerate in $p$, though this is not useful in our current work. Define the Lagrangian by the Legendre transform
\begin{equation*}
 L^c(x,t,\xi)
 :=\max_{0\le p\le C}\{\xi(p-c)-H(x,t,p)\}.
\end{equation*}
Then, as shown in Section \ref{sec:action}, for every admissible curve $\gamma(\cdot)$,
\begin{equation*}
 \sup_{0\le p(\cdot)\le C}\mathcal A^c(\gamma(\cdot),p(\cdot))
 =\int_0^tL^c(\gamma(s),s,\dot\gamma(s))\,\mathrm ds.
\end{equation*}
Hence the Hamiltonian  $\inf$-$\sup$ problem can be written exactly as a minimization problem for the Lagrangian action.  However, when $H$ is degenerate in $p$, the classical smooth Legendre correspondence between velocity and momentum may fail, and $L^c$ need not be differentiable.  The usual smooth Lagrangian formalism and Tonelli's calculus of variations therefore need not apply to the regularity problem for optimal trajectories.  We thus keep the momentum variable and work directly with the Hamiltonian  $\inf$-$\sup$ problem.

\subsection{Two difficulties in a direct method for the Hamiltonian action}\label{subsec:intro-obstructions}

There are two difficulties in applying a direct method to the Hamiltonian action.  First, a bounded optimizing sequence may give only
\[
 \gamma_n(\cdot)\to\gamma(\cdot)\quad\text{uniformly},\qquad
 \dot\gamma_n(\cdot)\stackrel{*}{\rightharpoonup}\dot\gamma(\cdot),
 \qquad
 p_n(\cdot)\stackrel{*}{\rightharpoonup}p(\cdot)
 \qquad\text{as }n\to\infty.
\]
These convergences are not enough to pass to the product $\dot\gamma_n(\cdot)p_n(\cdot)$ or to the nonlinear term $H(\gamma_n(\cdot),\cdot,p_n(\cdot))$ in the action. 
Second, the momentum constraint $0\le p\le C$ is part of the problem.  If the momentum component of an optimal trajectory meets $p=0$ or $p=C$, admissible variations may be one-sided.  Hence the usual two-sided first-variation argument does not directly give the second Hamilton equation.  This boundary difficulty has to be resolved before the regularity of an optimal trajectory is obtained.

\subsection{Main results and idea of proof}\label{subsec:intro-results}

As noted above, for a fixed terminal point $(t,x)$, the velocity constraint \eqref{eq:intro-admissible} restricts the possible initial points $y$.  Section~\ref{sec:reachable} shows that the corresponding reachable set $R_{t,x}$ is a compact interval. In the main theorem, Theorem \ref{thm:fixed-optimal},  we prove that,  for every $y\in R_{t,x}$, the $\inf$-$\sup$ problem \eqref{eq:intro-infsup} admits an optimal trajectory $(\gamma^*(\cdot),p^*(\cdot))\in C^2([0,t];\R\times[0,C])$ satisfying
\begin{equation*}
 \dot\gamma^*(s)=\partial_p H(\gamma^*(s),s,p^*(s)),\qquad
 \dot p^*(s)=-\partial_x H(\gamma^*(s),s,p^*(s)).
\end{equation*}
The theorem does not assert that every optimal trajectory has this regularity. 
We also prove the corresponding result for the Bolza problem: see Theorem \ref{thm:Bolza}.

The proof uses a finite-dimensional approximation together with a small quadratic regularization in the momentum variable, which makes the regularized Hamiltonian strictly convex in $p$. 
As in K.~Soga's finite-dimensional approximation method for Tonelli's calculus of variations \cite{Soga2023}, we first construct discrete optimizers and study their variational properties at the discrete level.  Once the discrete form of Hamilton's equations is obtained, these equations give uniform estimates for both the configuration and momentum components and hence the strong compactness needed to pass to the continuous limit.  In this way we avoid passing directly to the limit under the weak convergences described in the previous subsection. 

The finite-dimensional approximation also plays an essential role in the treatment of the momentum boundary.  One might try to handle the boundary directly by one-sided variations in continuous time. 
However, the set of times at which the momentum component of an optimizing trajectory reaches $0$ or $C$ may have a complicated structure: it may contain intervals or isolated contact times, and the momentum may enter and leave the boundary repeatedly.  
The finite-dimensional approximation avoids this difficulty by reducing the boundary analysis to finitely many momentum components.  
One-sided variations can then be used to control the boundary components and derive the discrete form of Hamilton's equations.

The interpolated discrete optimal trajectories converge to a classical trajectory satisfying Hamilton's equations. 
An additional comparison argument for an arbitrary admissible trajectory shows that the limiting Hamiltonian trajectory is optimal for the original  $\inf$-$\sup$ problem.

\subsection{Methodological perspectives and relation to known results}
\label{subsec:intro-literature}

Classical Tonelli's calculus of variations treats smooth Lagrangians defined on the whole velocity space and strictly convex and superlinear in the velocity variable.
The direct method gives minimizing curves in the class of absolutely continuous curves, and a regularity argument shows that these minimizers satisfy the Euler--Lagrange equation.
Through the Legendre transform they become classical Hamiltonian trajectories. 
At the same time, the value function associated with the Lagrangian variational problem solves the Hamilton-Jacobi equation generated by the corresponding Hamiltonian. 
This connection among variational minimizers, value functions, Hamilton--Jacobi equations, and Hamiltonian/Lagrangian dynamics is one of the foundations of weak KAM theory; see, for example, A.~Fathi \cite{Fathi2012}.

K.~Soga \cite{Soga2023} gave a finite-dimensional approximation method for Tonelli's calculus of variations that proves the existence and regularity of minimizers directly within the smooth class.
The finite-dimensional minimizers satisfy a discrete Euler-Lagrange equation, equivalently the discrete form of Hamilton's equations, and the corresponding polygons converge to smooth minimizing trajectories.
An important feature of this approach is that the variational equations and the resulting regularity estimates are obtained before the continuous limit is taken.
The present work follows this general idea, although the variational problem and the class of Lagrangians/Hamiltonians are quite different.

There are also variational approaches that start from Hamiltonian dynamics rather than from minimizers of the Lagrangian action.
A classical short-time result, going back to Weierstrass, states that, under suitable regularity and convexity assumptions, the projections of sufficiently short Hamiltonian trajectories onto configuration space minimize the Lagrangian action between their endpoints; see also A.~Fathi's discussion of small extremal curves \cite{Fathi2008}.
P.~Bernard \cite{Bernard2012} developed a Hamiltonian approach to the Lax--Oleinik operators and weak KAM solutions in the standard Tonelli setting, without relying on the Lagrangian.
In this setting, uniform convexity in the momentum variable yields the short-time twist property, so that the two configuration-space endpoints uniquely determine a Hamiltonian trajectory on a sufficiently short time interval.
The Hamiltonian action along this trajectory, regarded as a function of the two configuration-space endpoints, provides a generating function for the short-time Hamiltonian flow.
Its value coincides with the minimal Lagrangian action between these endpoints.

To extend this construction to a longer time interval, Bernard divides the interval into sufficiently short pieces and sums the corresponding generating functions.
The resulting finite-dimensional function depends on the intermediate configuration points.
Minimizing this function yields the minimal Lagrangian action between the prescribed endpoints over the whole time interval.
At a critical point of this function, the momentum at the end of each short trajectory equals the momentum at the start of the next one. 
Hence concatenating these short trajectories yields a single Hamiltonian trajectory over the whole time interval.
If the critical point is a minimum, the projection of this trajectory minimizes the Lagrangian action between the two endpoints. 
This construction is closely related to the broken-geodesics method introduced by M.~Chaperon \cite{Chaperon1984} and subsequently developed in the generating-function theory of J.-C.~Sikorav and C.~Viterbo \cite{Sikorav1987,Viterbo1992}.

P.~H.~Rabinowitz \cite{Rabinowitz1978} and F.~H.~Clarke and I.~Ekeland \cite{Clarke1979,ClarkeEkeland1980} developed another approach to Hamiltonian variational problems based on critical-point and minimax arguments.
In the Rabinowitz-type approach, the Hamiltonian action is a nonconvex, typically indefinite functional on phase space, and the configuration and momentum variables are varied independently.
Its critical points give classical Hamiltonian trajectories.
The dual-action methods of F.~H.~Clarke and I.~Ekeland provide related variational formulations, and related ideas were also extended to Hamiltonian inclusions \cite{Clarke1981}. 
A different variational framework based on selfdual Lagrangians
and convex duality was developed by N.~Ghoussoub and
A.~Moameni~\cite{GhoussoubMoameni2007,GhoussoubMoameni2008}
and N.~Ghoussoub~\cite{Ghoussoub2008}.
This approach obtains solutions of certain Hamiltonian systems
by minimizing suitable nonnegative functionals, rather than
through critical-point arguments for the classical Hamiltonian
action. 

More general nonsmooth theories of the calculus of variations take a different route.
F.~H.~Clarke's generalized Bolza problem and the convex theory of R.~T.~Rockafellar allow extended-valued Lagrangians (which may take the value $+\infty$) and constraints \cite{Clarke1976,Rockafellar1970,Rockafellar1971,Rockafellar1972}.
P.~D.~Loewen and R.~T.~Rockafellar \cite{LoewenRockafellar1996} obtained subgradient Euler--Lagrange and Hamiltonian conditions under very general assumptions; see also A.~D.~Ioffe and R.~T.~Rockafellar \cite{IoffeRockafellar1996} and the recent survey of A.~D.~Ioffe \cite{Ioffe2026}.
In these general settings, the momentum and Hamiltonian are treated in a generalized sense, and additional regularity may be needed to recover the classical Hamiltonian/Lagrangian dynamics.

To make our purposes clear, it is useful to distinguish three aspects of weak KAM theory.
The first is the classical dynamical aspect, in which minimizing or calibrated curves are represented by trajectories satisfying Hamilton's equations.
The second is the variational aspect, which studies minimizing, calibrated, static, and semistatic curves, without requiring them to be represented by classical Hamiltonian trajectories.

The third is the PDE/metric aspect, involving value functions, viscosity solutions, Lax--Oleinik operators, critical subsolutions, and projected Aubry sets.
Under the standard Tonelli assumptions, these three aspects are closely connected.
Under weaker regularity assumptions, however, the variational and PDE/metric aspects can still be studied even when a classical Hamiltonian flow is not available.
A.~Davini and M.~Zavidovique \cite{DaviniZavidovique2013} studied pairs of continuous convex Hamiltonians using Lax--Oleinik semigroups and investigated their weak KAM solutions and Aubry sets.
A.~Fathi and A.~Siconolfi \cite{FathiSiconolfi2005} developed a PDE approach to Aubry--Mather theory for continuous, coercive, and quasiconvex Hamiltonians without assuming a classical Hamiltonian flow.
A.~Agrachev and P.~W.~Y.~Lee \cite{AgrachevLee2010} also established a version of weak KAM theory for costs arising from affine control systems.

The study of the dynamical aspect of weak KAM theory goes beyond finding optimal values and optimal curves in the calculus of variations.
For instance, it also concerns the relation between the spatial derivatives of value functions and the momenta of optimal trajectories.
Since the Hamiltonian in our setting is smooth while the corresponding Lagrangian may be nonsmooth, it is natural and convenient to work directly on the Hamiltonian side.
With this motivation, the present paper provides a new direct method for the Hamiltonian $\inf$-$\sup$ action.
This method is intended as a starting point for developing the classical dynamical aspect of weak KAM theory in this degenerate setting. 
\medskip\medskip

\noindent {\bf Outline.} Section~\ref{sec:model-motivation} presents a concrete physical model leading to a Hamiltonian satisfying (H1)--(H3).
Section~\ref{sec:action} proves the existence of maximizing momentum curves for the Hamiltonian action.
Section~\ref{sec:reachable} characterizes the reachable set determined by the velocity constraint and proves the uniqueness of admissible curves at its boundary.
Section~\ref{sec:discrete} constructs the regularized finite-dimensional $\inf$-$\sup$ problem and proves the existence and regularity of its optimizers, deriving the discrete  Hamilton's equations.
Section~\ref{sec:optimality} passes to the limit in the discrete Hamilton's equations and proves the existence of an optimal  Hamiltonian trajectory for every reachable pair of endpoints.
Section~\ref{sec:bolza} extends this result to the Bolza $\inf$-$\sup$ problem.

\section{A concrete physical model}\label{sec:model-motivation}

The model comes from a reduction of the two-particle $\beta$-FPU oscillator to a nonautonomous system with one degree of freedom. 
The original two-degree-of-freedom Hamiltonian in the action-angle variables is given by  
\begin{equation*}
 \mathscr H(x,p)
 :=p_1+\omega p_2
 +\eps\bigl(p_1\sin^2x_1+\omega p_2\sin^2x_2\bigr)^2,
 \qquad \omega=\sqrt3,
\end{equation*}
for $x=(x_1,x_2)\in\R^2$ and $p=(p_1,p_2)\in[0,\infty)^2$, with $\eps>0$.  
Note that the momentum boundaries $p_1=0$ and $p_2=0$ are invariant under the Hamiltonian flow of $ \mathscr H$, since $\partial_{x_i}\mathscr H=0$ whenever $p_i=0$ ($i=1,2$).
Fix $h>0$ and consider the energy level $\{(x_1,x_2,p_1,p_2)\,|\,\mathscr H(x_1,x_2,p_1,p_2)=h\}$.  
Since
\begin{equation}\label{eq:model-p2-derivative}
 \partial_{p_2}\mathscr H(x_1,x_2,p_1,p_2)
 =\omega+2\eps\omega
   \bigl(p_1\sin^2x_1+\omega p_2\sin^2x_2\bigr)\sin^2x_2
 \ge \omega>0,
\end{equation}
the equation $\mathscr H(x,p)=h$ uniquely determines a nonnegative solution $p_2=p_2(x_1,x_2,p_1)$, which depends smoothly on $(x_1,x_2,p_1)$ in the domain  $x_1,x_2\in\R$ and $ 0\le p_1\le U(x_1;h)$,  where\footnote{On the energy surface $\mathscr H(x,p)=h$, the boundary $p_2=0$ corresponds to $p_1= U(x_1;h)$.} 
\begin{equation*}
 U(r;h):=\frac{2h}{1+\sqrt{1+4\eps h\sin^4r}}.
\end{equation*}
We set 
$$f(r,t,y):=-p_2(r,t,y),$$
which will be the reduced Hamiltonian, 
and rationalize the numerator to obtain\footnote{$y= U(r;h)$ is the positive solution of $h-(y+\eps y^2\sin^4 r)=0$.} 
\begin{equation}\label{eq:explicit-f}
\begin{aligned}
 f(r,t,y)
 ={}&-2\{h-(y+\eps y^2\sin^4 r)\}\\
 &\times\Bigl[
 \sqrt{\bigl(\omega+2\eps\omega y\sin^2r\sin^2t\bigr)^2
 +4\eps\omega^2\sin^4t\,\{h-(y+\eps y^2\sin^4r)\}}\\
 &\hspace{2.8cm}
 +\omega+2\eps\omega y\sin^2r\sin^2t
 \Bigr]^{-1},
 \qquad 0\le y\le U(r;h).
\end{aligned}
\end{equation}
The bracket in \eqref{eq:explicit-f} is at least $2\omega>0$ so that the formula is valid also when $\sin t=0$. Moreover,  $f(r,t,y)\le0$ throughout its domain.  Equivalently, $f$ is characterized by the energy identity
\begin{equation}\label{eq:implicit-f}
 \mathscr H\bigl(r,t,y,-f(r,t,y)\bigr)=h,
\end{equation}
or, more explicitly, 
\begin{equation}\label{eq:f-energy-identity}
 h=y-\omega f(r,t,y)
   +\eps\bigl(y\sin^2r-\omega f(r,t,y)\sin^2t\bigr)^2.
\end{equation}
The derivatives of $f$ can be easily obtained by implicit differentiation of  \eqref{eq:f-energy-identity}.   
Since $\dot x_2(s)=\partial_{p_2}\mathscr H(x_1(s),x_2(s),p_1(s),p_2(s))\ge\omega$, the map $s\mapsto x_2(s)$ is strictly increasing along every original Hamiltonian trajectory on $\{\mathscr H=h\}$.  Let $s=\eta(t)$ denote its inverse and set
\[
 r(t):=x_1(\eta(t)),\qquad y(t):=p_1(\eta(t)).
\]
We see that $( r(\cdot),y(\cdot))$ satisfies
\begin{equation}\label{eq:reduced-Ham}
 r'(t)=\partial_y f(r(t),t,y(t)),\qquad
 y'(t)=-\partial_r f(r(t),t,y(t)).
\end{equation}
In fact, the original Hamiltonian flow on the energy surface $\{\mathscr H=h\}$, with time variable $s$, is equivalent to the flow of \eqref{eq:reduced-Ham}. 
The equivalence is proved in Appendix~\ref{app:reduction}.

The reduced Hamiltonian has the following convexity property: 
\begin{equation*}
 \partial_y^2 f(r,t,y)
 =\frac{2\eps\bigl(\sin^2r-\omega\partial_y f(r,t,y)\sin^2t\bigr)^2}
 {\omega+2\eps\omega\bigl(y\sin^2r-\omega f(r,t,y)\sin^2t\bigr)\sin^2t}
 \ge0.
\end{equation*}
At $y=U(r;h)$, which corresponds to $p_2=0$ on the energy level $\{\mathscr H=h\}$, we have $f(r,t,U(r;h))=0$.  At $y=0$, which corresponds to $p_1=0$, \eqref{eq:explicit-f} shows that $f(r,t,0)$ is independent of $r$.  We denote this value by $f_0(t;h)$.

We now make a symplectic change of variables from $(r,y)$ to $(x,p)$ so that the variable momentum strip $0\le y\le U(r;h)$ becomes the fixed strip $0\le p\le C$ for some constant $C>0$.  
%This avoids a momentum constraint that depends on the position variable.
Let $C>0$ and let $a(\cdot)$ be a strictly increasing function.  Consider the change of variables
\begin{equation}\label{eq:sympl-map}
 \Phi(x,p)=(r,y),\qquad
 r:=a(x),\qquad y:=C^{-1}U(a(x);h)p,
 \qquad 0\le p\le C.
\end{equation}
For every $x$, this maps the fixed interval $0\le p\le C$ exactly onto $0\le y\le U(a(x);h)$.  We now choose $a(\cdot)$ and $C$ so that the change of variables is symplectic.  With
\[
 J:=\begin{pmatrix}0&1\\-1&0\end{pmatrix},
\qquad
 D\Phi(x,p)=
 \begin{pmatrix}
 a'(x)&0\\
 p\,\dfrac{\mathrm d}{\mathrm dx}\!\left(C^{-1}U(a(x);h)\right)&C^{-1}U(a(x);h)
 \end{pmatrix},
\]
a direct multiplication gives
\begin{equation}\label{eq:matrix-sympl}
 D\Phi(x,p)^T J D\Phi(x,p)
 =C^{-1}U(a(x);h)a'(x)\,J.
\end{equation}
Hence the symplectic condition is exactly
\begin{equation}\label{eq:sympl-cond}
 C^{-1}U(a(x);h)a'(x)=1.
\end{equation}
Set
\begin{equation*}
 C:=\frac1{2\pi}\int_0^{2\pi}U(r;h)\,\mathrm dr,
 \qquad
 V(r):=\int_0^r C^{-1}U(\rho;h)\,\mathrm d\rho,
 \qquad a:=V^{-1}.
\end{equation*}
Since $U(r;h)>0$, the function $V$ is strictly increasing.  Moreover, $V(r+2\pi)=V(r)+2\pi$, so $V$ maps $\R$ onto $\R$ and has a smooth inverse satisfying $a(x+2\pi)=a(x)+2\pi$.  Differentiating $V(a(x))=x$ yields \eqref{eq:sympl-cond}, and therefore \eqref{eq:matrix-sympl} shows that \eqref{eq:sympl-map} is symplectic.

Define
\begin{equation*}
 H(x,t,p)
 :=f\!\left(a(x),t,C^{-1}U(a(x);h)p\right).
\end{equation*}
Then $H\in C^2(\T\times\T\times[0,C];\R)$ is periodic in $(x,t)$ and
\begin{equation*}
 \partial_p^2 H(x,t,p)
 =\left(C^{-1}U(a(x);h)\right)^2
 \partial_y^2 f\!\left(a(x),t,C^{-1}U(a(x);h)p\right)\ge0.
\end{equation*}
Moreover, $H(x,t,C)=0$ and $H(x,t,0)=f_0(t;h)$.  Thus $H$  satisfies (H1)--(H3).

\section{Maximizing momentum curve}\label{sec:action}

Fix $c\in(0,C)$ and a terminal time $t>0$.  
A maximizing momentum curve associated with a fixed admissible curve $\gamma(\cdot)\in\Gamma_t$ exists but may fail to be unique.  
%Whenever pointwise values of $\dot\gamma$ are used in this subsection, we fix a measurable representative and redefine it on a null set, if necessary, so that the admissibility bounds in \eqref{eq:intro-admissible} hold for every $s\in[0,t]$.  This convention does not change the action.  The following result is the only measurable-selection result needed in the paper.
\begin{proposition}[Maximizing momentum curve]\label{prop:measurable-max}
For every $\gamma(\cdot)\in\Gamma_t$ there exists  $p_\gamma(\cdot)\in P_t$ such that
\begin{equation}\label{eq:sup-attained}
\Ac(\gamma(\cdot),p_\gamma(\cdot))
=\sup_{p(\cdot)\in P_t}\Ac(\gamma(\cdot),p(\cdot))=\int_0^t\max_{0\le p\le C}
\{\dot\gamma(s)(p-c)-H(\gamma(s),s,p)\}\dd s.
\end{equation}
For instance, with $E_\gamma:=\{s\in[0,t]\,|\,\partial_p H(\gamma(s),s,0)=\partial_p H(\gamma(s),s,C)\}$, one may choose $p_\gamma(s)=c$ for $s\in E_\gamma$, and for $s\notin E_\gamma$ one may take the left generalized inverse
\begin{equation}\label{eq:left-inverse}
p_\gamma(s)
:=\inf\{q\in[0,C]:\partial_p H(\gamma(s),s,q)\ge\dot\gamma(s)\}.
\end{equation}
%No condition on the size or shape of $E_\gamma$ is needed.
\end{proposition}

\begin{proof}
The set $E_\gamma$ is measurable.  Since  $p\mapsto \partial_p H(\gamma(s),s,p)$ is nondecreasing, $\partial_p H(\gamma(s),s,\cdot)$ is constant on $[0,C]$ for each $s\in E_\gamma$. 
 Hence the admissibility implies that $\dot\gamma(s)=\partial_p H(\gamma(s),s,0)$ and $p\mapsto \dot\gamma(s)(p-c)-H(\gamma(s),s,p)$ is constant on $[0,C]$.   
Therefore every $p\in[0,C]$ is a maximizer for $ \dot\gamma(s)(\cdot-c)-H(\gamma(s),s,\cdot)$.  We choose $p_\gamma(s):=c$ for each $s\in E_\gamma$.

For almost every $s\in[0,t]\setminus E_\gamma$, define $p_\gamma(s)$ by \eqref{eq:left-inverse}. The admissibility and continuity of $\partial_p H$ yield $\partial_p H(\gamma(s),s,p_\gamma(s))=\dot\gamma(s)$. By the convexity of $H$ in $p$, this shows that $p_\gamma(s)$ maximizes $\dot\gamma(s)(p-c)-H(\gamma(s),s,p)$ over $p\in[0,C]$. 
For $r\in[0,C)$, the nondecreasing property of $\partial_pH(\gamma(s),s,\cdot)$ implies that 
\[
\{s\in [0,t]\setminus E_\gamma\,|\,p_\gamma(s)>r\}
=\{s\in [0,t]\setminus E_\gamma\,|\,\partial_p H(\gamma(s),s,r)<\dot\gamma(s)\},
\]
which is measurable.  
Thus $p_\gamma$ is measurable on $[0,t]\setminus E_\gamma$. Since $p_\gamma=c$ on the measurable set $E_\gamma$, $p_\gamma(\cdot)$ is measurable on $[0,t]$.

By the preceding construction, for almost every $s\in[0,t]$,
\[
\max_{0\le p\le C}
\{\dot\gamma(s)(p-c)-H(\gamma(s),s,p)\}
=
\dot\gamma(s)(p_\gamma(s)-c)
-H(\gamma(s),s,p_\gamma(s)),
\]
and the right-hand side is measurable.
For any $p(\cdot)\in P_t$ and for almost every $s\in[0,t]$, we have
\[
\dot\gamma(s)(p(s)-c)-H(\gamma(s),s,p(s))
\le
\dot\gamma(s)(p_\gamma(s)-c)
-H(\gamma(s),s,p_\gamma(s)).
\]
Integration gives \eqref{eq:sup-attained}.
\end{proof}
%\begin{remark}[Partial degeneracy]\label{rem:partial-deg}
%Suppose that, for some $(x,s,\xi)$, the set $\{p\in[0,C]:\partial_pH(x,s,p)=\xi\}$ is the interval $[p_-,p_+]$, where $p_-<p_+$.  Then every $p\in[p_-,p_+]$ maximizes $p\mapsto \xi(p-c)-H(x,s,p)$.
%In this case the choice in Proposition~\ref{prop:measurable-max} gives the left endpoint $p_-$.  Thus degeneracy may destroy uniqueness, but a measurable maximizing momentum curve still exists for a fixed admissible curve $\gamma(\cdot)$.
%\end{remark}

\section{ The reachable set and endpoint uniqueness}\label{sec:reachable}

Fix $x\in\R$ and $t>0$, and impose the terminal condition $\gamma(t)=x$.  In the two-fixed-endpoint problem, the velocity constraint \eqref{eq:intro-admissible} does not allow an arbitrary initial point $y$.  We therefore define the reachable set by
\[
R_{t,x}
:=\left\{y\in\R\,|\,\text{there exists }\gamma(\cdot)\in\Gamma_t\text{ with }\gamma(0)=y\text{ and }\gamma(t)=x\right\}.
\]
The two boundary velocities in \eqref{eq:intro-admissible} give the backward terminal-value problems
\begin{align}
\dot{\overline{\gamma}}(s)&=\partial_p H(\overline{\gamma}(s),s,C),&\overline{\gamma}(t)&=x,\label{eq:extreme-C}\\
\dot{\underline{\gamma}}(s)&=\partial_p H(\underline{\gamma}(s),s,0),&\underline{\gamma}(t)&=x.\label{eq:extreme-0}
\end{align}
Backward comparison gives
\[
\overline{\gamma}(0)\le\underline{\gamma}(0).
\]

\begin{proposition}\label{prop:reachable}
One has
\[
R_{t,x}=[\overline{\gamma}(0),\underline{\gamma}(0)].
\]
Moreover, if $\gamma(\cdot)\in\Gamma_t$ satisfies $\gamma(t)=x$ and $\gamma(0)=\overline{\gamma}(0)$, then $\gamma(s)=\overline{\gamma}(s)$ for all $s\in[0,t]$.  If instead $\gamma(0)=\underline{\gamma}(0)$, then $\gamma(s)=\underline{\gamma}(s)$ for all $s\in[0,t]$.
\end{proposition}

\begin{proof}
Let $y\in R_{t,x}$ and choose $\gamma(\cdot)\in\Gamma_t$ with $\gamma(0)=y$ and $\gamma(t)=x$.  
Since $\dot\gamma(s)\le \partial_p H(\gamma(s),s,C)$, we have 
$$0\ge \dot\gamma(s)-\partial_p H(\gamma(s),s,C)  
=  \dot\gamma(s) - \dot{\overline{\gamma}}(s)+\partial_p H(\overline{\gamma}(s),s,C)- \partial_p H(\gamma(s),s,C)  $$
The standard backward Gronwall inequality gives $\gamma(s)\ge\overline{\gamma}(s)$.  Similarly, $\dot\gamma(s)\ge \partial_p H(\gamma(s),s,0)$ implies $\gamma(s)\le\underline{\gamma}(s)$.  Thus $y=\gamma(0)\in[\overline{\gamma}(0),\underline{\gamma}(0)]$, and hence $R_{t,x}\subset[\overline{\gamma}(0),\underline{\gamma}(0)]$.

Conversely, for each constant $p\in[0,C]$ solve
\begin{equation*}
\dot\gamma(s;p)=\partial_p H(\gamma(s;p),s,p),
\qquad \gamma(t;p)=x. 
\end{equation*}
By (H2), each $\gamma(\cdot;p)$ belongs to $\Gamma_t$. 
The map $p\mapsto\gamma(0;p)$ is continuous, with  $\gamma(0;0)=\underline{\gamma}(0)$ and $\gamma(0;C)=\overline{\gamma}(0)$.   Therefore the  intermediate value theorem implies that for every $y\in[\overline{\gamma}(0),\underline{\gamma}(0)]$ there exists $p\in[0,C]$ such that $\gamma(0;p)=y$.  This proves $R_{t,x}\supset [\overline{\gamma}(0),\underline{\gamma}(0)]$.

Suppose $\gamma(\cdot)\in\Gamma_t$ satisfies $\gamma(t)=x$ and $\gamma(0)=\overline{\gamma}(0)$.  Put $z(s):=\gamma(s)-\overline{\gamma}(s)$ and $\delta(s):=\partial_p H(\gamma(s),s,C)-\dot\gamma(s)=    \dot{\overline{\gamma}}(s) -\dot\gamma(s)+ \partial_p H(\gamma(s),s,C) - \partial_p H(\overline{\gamma}(s),s,C)  \ge0$.
By the mean value theorem, we find  $\alpha(\cdot)\in L^\infty(0,t)$ such that 
\[
\partial_p H(\gamma(s),s,C)-\partial_p H(\overline{\gamma}(s),s,C)=\alpha(s)z(s).
\]
Hence we have $z'(s)=\alpha(s)z(s)-\delta(s)$ with $z(t)=0$, which leads to  
\[
z(s)=\int_s^t
 \exp\left(-\int_s^r\alpha(\sigma)\dd\sigma\right)\delta(r)\dd r\ge0.
\]
Since also $z(0)=0$, one obtains $\delta(s)=0$ for almost every $s$, and hence $\dot\gamma(s)=\partial_p H(\gamma(s),s,C)$ for almost every $s$.  Uniqueness of the terminal-value problem gives $\gamma(s)=\overline{\gamma}(s)$ for every $s$.  The case $\gamma(0)=\underline{\gamma}(0)$ is analogous.
\end{proof}

\section{Finite-dimensional approximation of the  inf-sup problem}\label{sec:discrete}

Fix $c\in(0,C)$ and a terminal time $t>0$.  We construct a finite-dimensional approximation of the two-fixed-endpoint  $\inf$-$\sup$ problem.  We first introduce a regularized discrete problem and characterize its reachable set.  We then obtain discrete optimizers, study one-sided discrete variations with a careful treatment of  the momentum boundary, and derive the discrete form of Hamilton's equations for these optimizers. 

\subsection{The regularized discrete problem}\label{dis-pro}

For each $K\in\mathbb N$,  define 
\begin{align*}
&\tau:=\frac{t}{K},\quad s_k:=k\tau,\\
&\gamma=(\gamma_0,\ldots,\gamma_K)\in\R^{K+1},\quad p=(p_1,\ldots,p_K)\in [0,C]^K,\quad
\dot{\gamma}_k:=\frac{\gamma_k-\gamma_{k-1}}{\tau}.
\end{align*}
We regularize the Hamiltonian by $H_\tau(x,s,p):=H(x,s,p)+\frac{\sqrt\tau}{2}(p-c)^2$.  Then
\begin{equation}\label{eq:Htau-strict}
\partial_p^2 H_\tau=\partial_p^2 H+\sqrt\tau>0.
\end{equation}
The regularization has two purposes. It makes the maximizing momentum unique at the discrete level, and it slightly enlarges the admissible discrete velocity interval, which will be used in Section \ref{subsec:optimality-limit} to approximate continuous admissible curves by discrete admissible curves.
Note that the proof does not require the particular factor $\sqrt\tau$.  It works with any positive number $r_\tau$ such that $r_\tau\to0$ and $\tau/r_\tau\to0$ as $\tau\to0+$. 

The discrete Hamiltonian action is defined as 
\begin{equation*}
\AcK(\gamma,p)
:=\sum_{k=1}^K
\left\{\dot{\gamma}_k(p_k-c)-H(\gamma_k,s_k,p_k)
-\frac{\sqrt\tau}{2}(p_k-c)^2\right\}\tau.
\end{equation*}
The discrete admissibility condition is
\begin{equation}\label{eq:disc-adm}
 \partial_p H(\gamma_k,s_k,0)-\sqrt\tau\,c
 \le \dot{\gamma}_k\le
 \partial_p H(\gamma_k,s_k,C)+\sqrt\tau\,(C-c),
 \qquad k=1,\ldots,K.
\end{equation}
For every discrete admissible curve $\gamma$, strict monotonicity in \eqref{eq:Htau-strict} gives a unique $p_k^\gamma\in[0,C]$ such that
\begin{equation}\label{eq:disc-velocity}
 \dot{\gamma}_k=\partial_p H(\gamma_k,s_k,p_k^\gamma)+\sqrt\tau\,(p_k^\gamma-c).
\end{equation}
Equivalently, $p_k^\gamma$ uniquely maximizes $p\mapsto \dot{\gamma}_k(p-c)-H(\gamma_k,s_k,p)-\frac{\sqrt\tau}{2}(p-c)^2$ on $[0,C]$.  Set $p^\gamma:=(p_1^\gamma,\ldots,p_K^\gamma)$; then $\max_{p\in[0,C]^K}\AcK(\gamma,p)=\AcK(\gamma,p^\gamma)$.
The map $\gamma\mapsto p^\gamma$ is continuous on the set of discrete admissible curves, regarded as a subset of $\mathbb R^{K+1}$ with the Euclidean topology.  Hence the maximized discrete action is also continuous there.

\subsection{Discrete reachable set}

For the regularized discrete problem, define
\[
R_{t,x}^K
:=\left\{y\in\R\,|\,\text{there exists a discrete admissible curve }\gamma\text{ with }\gamma_0=y\text{ and }\gamma_K=x\right\}.
\]
Define the regularized boundary backward sequences recursively by 
\begin{align}
\overline{\gamma}_{k-1}^{\,K}
&:=\overline{\gamma}_k^{\,K}
-\tau\{\partial_p H(\overline{\gamma}_k^{\,K},s_k,C)+\sqrt\tau(C-c)\}\mbox{\,\,\, for  $k=1,\ldots,K$ with $\overline{\gamma}_K^{\,K}:=x$},\label{eq:disc-ext-C}\\
\underline{\gamma}_{k-1}^{\,K}
&:=\underline{\gamma}_k^{\,K}
-\tau\{\partial_p H(\underline{\gamma}_k^{\,K},s_k,0)-\sqrt\tau c\}\mbox{\,\,\, for  $k=1,\ldots,K$ with $\underline{\gamma}_K^{\,K}:=x$}.\label{eq:disc-ext-0}
\end{align}

\begin{lemma}[The discrete reachable set]\label{lem:disc-reachable}
For all sufficiently large $K\in\mathbb N$,
\[
 R_{t,x}^K
 = [\overline{\gamma}_0^{\,K},\underline{\gamma}_0^{\,K}].
\]
Moreover,
\begin{equation}\label{eq:extreme-conv}
 \overline{\gamma}_0^{\,K}\to\overline{\gamma}(0),
 \qquad
 \underline{\gamma}_0^{\,K}\to\underline{\gamma}(0)
 \qquad\text{as }K\to\infty.
\end{equation}
In fact, for all sufficiently large $K$,
\begin{equation}\label{eq:reachable-strict-inclusion}
 \overline{\gamma}_0^{\,K}
 <
 \overline{\gamma}(0)
 \le
 \underline{\gamma}(0)
 <
 \underline{\gamma}_0^{\,K},
\end{equation}
and consequently $R_{t,x}\subset\operatorname{int}R_{t,x}^K$.
\end{lemma}

\begin{proof}
Set
\[
 B_k(z):=z-\tau\{\partial_p H(z,s_k,C)+\sqrt\tau(C-c)\},
 \qquad
 A_k(z):=z-\tau\{\partial_p H(z,s_k,0)-\sqrt\tau c\}.
\]
For all sufficiently large $K$, equivalently for sufficiently small $\tau$, both maps are strictly increasing because
\begin{align}\label{CFL}
 B_k'(z)=1-\tau \partial_{xp}^2 H(z,s_k,C)>0,
 \qquad
 A_k'(z)=1-\tau \partial_{xp}^2 H(z,s_k,0)>0.
\end{align}
If $\gamma$ is a discrete admissible curve with $\gamma_K=x$ and $\gamma_0=y$, then in view of \eqref{eq:disc-velocity} we have 
\[
 \gamma_{k-1}=\gamma_k-\tau \dot{\gamma}_k=\gamma_k-\tau \{\partial_p H(\gamma_k,s_k,p_k^\gamma)+\sqrt\tau\,(p_k^\gamma-c)\}
 \ge B_k(\gamma_k),
 \qquad
 \gamma_{k-1}\le A_k(\gamma_k).
\]
Since \eqref{eq:disc-ext-C} and \eqref{eq:disc-ext-0} are written as $\overline{\gamma}_{k-1}^{\,K}=B_k(\overline{\gamma}_{k}^{\,K})$,  $\underline{\gamma}_{k-1}^{\,K}=A_k(\underline{\gamma}_{k}^{\,K})$, backward induction together with \eqref{CFL} gives
$\overline{\gamma}_k^{\,K}\le\gamma_k\le\underline{\gamma}_k^{\,K}$ for every $k$.  Hence every $y\in R_{t,x}^K$ belongs to $[\overline{\gamma}_0^{\,K},\underline{\gamma}_0^{\,K}]$.

For $q\in[0,C]$ define
\begin{equation*}
 \gamma_{k-1}^{q,K}
 :=\gamma_k^{q,K}-\tau\{\partial_p H(\gamma_k^{q,K},s_k,q)+\sqrt\tau(q-c)\},
 \qquad \gamma_K^{q,K}:=x.
\end{equation*}
Then  $\gamma^{q,K}=( \gamma_{0}^{q,K},\ldots, \gamma_{K}^{q,K}) $ is a   discrete admissible curve. 
Its initial value $\gamma^{q,K}_0$  depends continuously on $q$ with  $\gamma^{C,K}_0=\overline{\gamma}_0^{\,K}$ and $\gamma^{0,K}_0=\underline{\gamma}_0^{\,K}$.  
The intermediate value theorem implies that  $[\overline{\gamma}_0^{\,K},\underline{\gamma}_0^{\,K}]\subset R_{t,x}^K$. Thus $R_{t,x}^K
 = [\overline{\gamma}_0^{\,K},\underline{\gamma}_0^{\,K}]$. 
 
 The classical result on convergence of the Euler-type method for ODEs shows that
\eqref{eq:disc-ext-C}, \eqref{eq:disc-ext-0} tend to
\eqref{eq:extreme-C}, \eqref{eq:extreme-0} as $K\to\infty$, i.e.,
\begin{equation*}
 \max_k
 |\overline{\gamma}_k^{\,K}-\overline{\gamma}(s_k)|
 =
 O(\tau+\sqrt\tau),
 \qquad
 \max_k
 |\underline{\gamma}_k^{\,K}-\underline{\gamma}(s_k)|
 =
 O(\tau+\sqrt\tau).
\end{equation*}
In particular, \eqref{eq:extreme-conv} holds. 
We refine this estimate.
Since $\overline{\gamma}$ solves \eqref{eq:extreme-C}, Taylor's formula
gives, uniformly in $k$,
\begin{equation*}
 \overline{\gamma}(s_{k-1})
 =\overline{\gamma}(s_k)
 - \tau\partial_pH(\overline{\gamma}(s_k),s_k,C)
 + O(\tau^2).
\end{equation*}
Combining this with \eqref{eq:disc-ext-C}, we obtain
\begin{align*}
 &\overline{\gamma}_{k-1}^{\,K}
 - \overline{\gamma}(s_{k-1})\\
 &\quad =
 \overline{\gamma}_k^{\,K}
 - \overline{\gamma}(s_k)
 -\tau\{ \partial_pH(\overline{\gamma}_k^{\,K},s_k,C)
 - \partial_pH(\overline{\gamma}(s_k),s_k,C)\}
  -\tau\sqrt\tau(C-c)
 + O(\tau^2).
\end{align*}
For all sufficiently small $\tau$, we obtain by backward iteration with  
$\overline{\gamma}_K^{\,K}=\overline{\gamma}(t)=x$, 
\begin{equation*}
 \overline{\gamma}_0^{\,K}
 - \overline{\gamma}(0)
 \le-e^{-2t\|\partial_{xp}^2H\|_{L^\infty}}
 t(C-c)\sqrt\tau
 + O(\tau).
\end{equation*}
Since $C-c>0$, it follows that $ \overline{\gamma}_0^{\,K}
 < \overline{\gamma}(0)$ for all sufficiently large $K$ (i.e., all sufficiently small $\tau$).
Similarly,  \eqref{eq:extreme-0} and  \eqref{eq:disc-ext-0}  give
\begin{align*}
 &\underline{\gamma}_{k-1}^{\,K}
 - \underline{\gamma}(s_{k-1})\\
 &\quad=
 \underline{\gamma}_k^{\,K}
 - \underline{\gamma}(s_k)
 - \tau\{
 \partial_pH(\underline{\gamma}_k^{\,K},s_k,0)
 - \partial_pH(\underline{\gamma}(s_k),s_k,0)\}
  + \tau\sqrt\tau c
 +O(\tau^2).
\end{align*}
The same backward iteration with
$\underline{\gamma}_K^{\,K}=\underline{\gamma}(t)=x$ yields
\begin{equation*}
 \underline{\gamma}_0^{\,K}
 - \underline{\gamma}(0)
 \ge e^{-2t\|\partial_{xp}^2H\|_{L^\infty}}
 tc\sqrt\tau
 + O(\tau).
\end{equation*}
Since $c>0$, we obtain $ \underline{\gamma}_0^{\,K}
 > \underline{\gamma}(0)$ for all sufficiently large $K$  (i.e., all sufficiently small $\tau$). 
Since $R_{t,x} = [\overline{\gamma}(0),\underline{\gamma}(0)]$,  we conclude 
$ R_{t,x} \subset\operatorname{int}R_{t,x}^K$.
\end{proof}

\subsection{Discrete optimizers}
Fix any $y\in R_{t,x}^K$.  Let $G_K(y,x)$ be the set of discrete admissible curves with $\gamma_0=y$ and $\gamma_K=x$.  The set $G_K(y,x)$ is nonempty, closed, and bounded, hence compact.  Since $\gamma\mapsto\max_{p\in[0,C]^K}\AcK(\gamma,p)=\AcK(\gamma,p^\gamma)$ is continuous, it has a minimizer $\gamma^{K,*}$ for 
\begin{align}\label{discrete-inf-sup}
\inf_{\substack{\gamma: \,\text{admissible} \\
                   \gamma_0=y,\ \gamma_K=x}}
\left\{ \max_{p\in[0,C]^K}\AcK(\gamma,p)\right\}=\inf_{\substack{\gamma: \,\text{admissible} \\
                   \gamma_0=y,\ \gamma_K=x}}
\AcK(\gamma,p^\gamma). 
\end{align} 
 We denote the associated unique maximizing momenta by $p_k^{K,*}:=p_k^{\gamma^{K,*}}$.
Thus $(\gamma^{K,*},p^{K,*})$ is a discrete optimizing trajectory for the left-hand side of \eqref{discrete-inf-sup}.

\subsection{One-sided directional derivatives}\label{subsec:one-sided}

Because the velocity constraints in \eqref{eq:disc-adm} may be fulfilled at the boundary, a discrete variation of a discrete admissible curve need not be admissible in both directions.  We therefore first derive a formula for the one-sided directional derivative of the maximized discrete action $\AcK(\gamma,p^\gamma)$. 

\begin{lemma}[One-sided derivative of the maximized discrete action]\label{lem:envelope}
Let $\gamma=(\gamma_0,\ldots,\gamma_K)$ be a discrete admissible curve and let $\omega=(\omega_0,\ldots,\omega_K)\in\R^{K+1}$.  Set $\gamma^\eta:=\gamma+\eta\omega$ for $\eta\in\R$.  
If $\gamma^\eta$ is admissible for all sufficiently small $\eta\ge0$, then
\[
 \left.\frac{\dd}{\dd\eta}
 \max_{p\in[0,C]^K}\AcK(\gamma^\eta,p)\right|_{\eta=0+}
 =\left.\frac{\dd}{\dd\eta}
 \AcK(\gamma^\eta,p^{\gamma^\eta})\right|_{\eta=0+}
 =\left.\frac{\dd}{\dd\eta}
 \AcK(\gamma^\eta,p^\gamma)\right|_{\eta=0+},
\]
where the second equality is the assertion of the lemma. 
%The first equality is an identity by the definition of $p^{\gamma^\eta}$.  The content of the lemma is the second equality: at $\eta=0$, the maximizing momentum vector may be kept fixed at $p^\gamma$ when taking the one-sided derivative.  
The analogous three-term equality with $0-$ holds if $\gamma^\eta$ is admissible for all sufficiently small $\eta\le0$.
\end{lemma}

\begin{proof}
Suppose that $\gamma^\eta$ is admissible for all sufficiently small $\eta\ge0$. 
For $\eta>0$ small, maximality of $p^{\gamma^\eta}$ for $\gamma^\eta$ and of $p^\gamma$ for $\gamma$ gives
\begin{align*}
 &\AcK(\gamma^\eta,p^\gamma)-\AcK(\gamma,p^\gamma)
 \le \AcK(\gamma^\eta,p^{\gamma^\eta})-\AcK(\gamma,p^\gamma)\\
 &\qquad=
 \max_{p\in[0,C]^K}\AcK(\gamma^\eta,p)
 -\max_{p\in[0,C]^K}\AcK(\gamma,p)\le
 \AcK(\gamma^\eta,p^{\gamma^\eta})-\AcK(\gamma,p^{\gamma^\eta}).
\end{align*}
Dividing by $\eta>0$, we obtain
\begin{align*}
 &\frac{\AcK(\gamma^\eta,p^\gamma)-\AcK(\gamma,p^\gamma)}{\eta}
 \le
 \frac{\displaystyle
 \max_{p\in[0,C]^K}\AcK(\gamma^\eta,p)
 -\max_{p\in[0,C]^K}\AcK(\gamma,p)}{\eta}
 \le\frac{\AcK(\gamma^\eta,p^{\gamma^\eta})
 -\AcK(\gamma,p^{\gamma^\eta})}{\eta}.
\end{align*}
We see that the two outer difference quotients have the same limit as $\eta\to0+$.  Indeed, since $\gamma_k^\eta=\gamma_k+\eta\omega_k$ and the corresponding discrete velocity is 
$\dot{\gamma}_k^\eta=\dot{\gamma}_k+\frac{\eta}{\tau}(\omega_k-\omega_{k-1})$, 
the definition of $\AcK$ yields,  
%\begin{align*}
% \frac{\AcK(\gamma^\eta,p^\gamma)-\AcK(\gamma,p^\gamma)}{\eta}
% =\sum_{k=1}^K(\omega_k-\omega_{k-1})(p_k^\gamma-c)
%-\tau\sum_{k=1}^K
% \frac{H(\gamma_k+\eta\omega_k,s_k,p_k^\gamma)
 %      -H(\gamma_k,s_k,p_k^\gamma)}{\eta}.
%\end{align*}
as $\eta\to0+$,
\begin{align*}
 \frac{\AcK(\gamma^\eta,p^\gamma)-\AcK(\gamma,p^\gamma)}{\eta}
 \longrightarrow
 \sum_{k=1}^K(\omega_k-\omega_{k-1})(p_k^\gamma-c)
 -\tau\sum_{k=1}^K
 \partial_xH(\gamma_k,s_k,p_k^\gamma)\omega_k.
\end{align*}
The continuity of $\gamma\mapsto p^\gamma$ gives $p_k^{\gamma^\eta}\to p_k^\gamma$ as $\eta\to0+$ for each $k$. 
Therefore, with the continuity of $\partial_xH$, we obtain,  
as $\eta\to0+$, 
%\begin{align*}
% &\frac{\AcK(\gamma^\eta,p^{\gamma^\eta})
% -\AcK(\gamma,p^{\gamma^\eta})}{\eta}\\
% &=\sum_{k=1}^K(\omega_k-\omega_{k-1})(p_k^{\gamma^\eta}-c)\\
 %&\quad-\tau\sum_{k=1}^K
% \frac{H(\gamma_k+\eta\omega_k,s_k,p_k^{\gamma^\eta})
%       -H(\gamma_k,s_k,p_k^{\gamma^\eta})}{\eta}.
%\end{align*}
%\[
% \frac{H(\gamma_k+\eta\omega_k,s_k,p_k^{\gamma^\eta})
%       -H(\gamma_k,s_k,p_k^{\gamma^\eta})}{\eta}
% \longrightarrow
% \partial_xH(\gamma_k,s_k,p_k^\gamma)\omega_k.
%\]
\begin{align*}
 \frac{\AcK(\gamma^\eta,p^{\gamma^\eta})
 -\AcK(\gamma,p^{\gamma^\eta})}{\eta}\longrightarrow
 \sum_{k=1}^K(\omega_k-\omega_{k-1})(p_k^\gamma-c)
 -\tau\sum_{k=1}^K
 \partial_xH(\gamma_k,s_k,p_k^\gamma)\omega_k.
\end{align*}
This proves the second equality and the existence of the right derivative.  A parallel argument deals with the case  $\eta\le0$.  
\end{proof}

\subsection{Boundary transitions}

We begin the proof that, for an interior reachable point $y\in {\rm int} R^K_{t,x}$, the maximizing momentum $p^{K,*}$ of every discrete optimizer stays in the interior of $[0,C]^K$. This interiority is crucial because it allows two-sided variations  along the minimizer $\gamma^{K,*}$ leading to discrete Hamilton's equations.  
Our first step is to rule out a direct transition between $0$ and $C$ for the discrete maximizing momentum.  
Using this, we then show that an interior reachable point gives at least one momentum component belonging to $(0,C)$, and finally propagate this interiority to all maximizing momentum components.

In the boundary arguments below, we take $K$ sufficiently large that
\begin{equation}\label{eq:tau-Hxp}
\tau\|\partial_{xp}^2 H\|_{L^\infty}<1.
\end{equation}
This condition guarantees that one-sided node variations\footnote{ For a discrete admissible curve  $\gamma=(\gamma_0,\ldots,\gamma_K)$ and $\omega=(0,\ldots,0,\omega_k,0,\ldots,0)\in\R^{K+1}$, the discrete curve $\gamma^\eta:=\gamma+\eta\omega$ for $\eta\in\R$ is called a $k$-node variation of $\gamma$.} used below move the velocity constraints achieved at the boundary strictly inward. 

\begin{lemma}[No direct transition between $0$ and $C$]\label{lem:no-switch}
Let $(\gamma^{K,*},p^{K,*})$ be a discrete optimizing trajectory for the left-hand side of \eqref{discrete-inf-sup}. Then the momentum curve $p^{K,*}$ cannot satisfy
$p_j^{K,*}=0$, $p_{j+1}^{K,*}=C$, nor
$p_j^{K,*}=C$, $p_{j+1}^{K,*}=0$,
for any $j=1,\ldots,K-1$.
\end{lemma}

\begin{proof} Assume first that $p_j^{K,*}=0$ and $p_{j+1}^{K,*}=C$ for some $j\in\{1,\ldots,K-1\}$. 
We vary only the node $\gamma_j^{K,*}$:
\[
\gamma_j^\eta:=\gamma_j^{K,*}+\eta,
\qquad
\gamma_i^\eta:=\gamma_i^{K,*}\quad(i\ne j).
\]
Then
\[
\dot{\gamma}_j^\eta
=\dot{\gamma}_j^{K,*}+\frac{\eta}{\tau},
\qquad
\dot{\gamma}_{j+1}^\eta
=\dot{\gamma}_{j+1}^{K,*}-\frac{\eta}{\tau}.
\]
For $\eta=0$, the $j$-th velocity lies at its lower regularized boundary as 
\[
\dot{\gamma}_j^\eta|_{\eta=0}=\dot{\gamma}_j^{K,*}
=
\partial_p H(\gamma_j^{K,*},s_j,0)-\sqrt\tau c.
\]
For $\eta>0$, \eqref{eq:tau-Hxp} gives 
\[
\frac{\dd}{\dd\eta}
\left[
\dot{\gamma}_j^\eta
-\partial_pH(\gamma_j^\eta,s_j,0)
+\sqrt\tau c
\right]_{\eta=0}
=
\frac1\tau-\partial_{xp}^2H(\gamma_j^{K,*},s_j,0)>0.
\]
Hence, for all sufficiently small $\eta>0$, we have 
$$\dot{\gamma}_j^\eta
-(\partial_pH(\gamma_j^\eta,s_j,0)
-\sqrt\tau c)>\Big[\dot{\gamma}_j^\eta
-(\partial_pH(\gamma_j^\eta,s_j,0)
-\sqrt\tau c)\Big]\Big|_{\eta=0}=0.$$
Hence, for all sufficiently small $\eta>0$, we have
\[
\dot{\gamma}_j^\eta
>
\partial_pH(\gamma_j^\eta,s_j,0)-\sqrt\tau c.
\]
Moreover, at $\eta=0$, we have by (H2), 
\[
\partial_pH(\gamma_j^{K,*},s_j,C)+\sqrt\tau(C-c)
-\dot{\gamma}_j^{K,*}
=
\partial_pH(\gamma_j^{K,*},s_j,C)
-\partial_pH(\gamma_j^{K,*},s_j,0)
+\sqrt\tau C
\ge\sqrt\tau C>0.
\]
Hence, by continuity, for all sufficiently small $\eta>0$,
\[
\dot{\gamma}_j^\eta<
\partial_pH(\gamma_j^\eta,s_j,C)+\sqrt\tau(C-c).
\]
At the same time, since $p_{j+1}^{K,*}=C$, the $(j+1)$-st velocity lies at
its upper regularized boundary,
\[
\dot{\gamma}_{j+1}^\eta|_{\eta=0}=\dot{\gamma}_{j+1}^{K,*}
=
\partial_pH(\gamma_{j+1}^{K,*},s_{j+1},C)
+\sqrt\tau(C-c).
\]
Since $\gamma_{j+1}^\eta=\gamma_{j+1}^{K,*}$ and
$\dot{\gamma}_{j+1}^\eta
=\dot{\gamma}_{j+1}^{K,*}-\eta/\tau$, we have
\[
\frac{\dd}{\dd\eta}
\left[
\partial_pH(\gamma_{j+1}^\eta,s_{j+1},C)
+\sqrt\tau(C-c)
-\dot{\gamma}_{j+1}^\eta
\right]_{\eta=0}
=
\frac1\tau>0.
\]
Hence, for all sufficiently small $\eta>0$, it holds that 
$$\partial_pH(\gamma_{j+1}^\eta,s_{j+1},C)
+\sqrt\tau(C-c)
-\dot{\gamma}_{j+1}^\eta>\Big[\partial_pH(\gamma_{j+1}^\eta,s_{j+1},C)
+\sqrt\tau(C-c)
-\dot{\gamma}_{j+1}^\eta\Big]\Big|_{\eta=0}=0,$$
and therefore, $\partial_pH(\gamma_{j+1}^\eta,s_{j+1},0)
-\sqrt\tau c < \dot{\gamma}_{j+1}^\eta<\partial_pH(\gamma_{j+1}^\eta,s_{j+1},C)
+\sqrt\tau(C-c)$. 
Since $\dot{\gamma}_{k}^\eta=\dot{\gamma}_{k}^{K,*}$ for all $k\neq j, j+1$, we confirm that $\gamma^\eta$ is admissible for all sufficiently small
$\eta>0$ with $\gamma^\eta_K=x$ and $\gamma^\eta_0=y$. 

We can therefore apply  Lemma~\ref{lem:envelope}: on either side for which the variation is admissible,
\begin{equation}\label{eq:first-var}
\left.\frac{\dd}{\dd\eta}
\max_{p\in[0,C]^K}\AcK(\gamma^\eta,p)\right|_{\eta=0\pm}
=\left.\frac{\dd}{\dd\eta}\AcK(\gamma^\eta,p^{\gamma^\eta})\right|_{\eta=0\pm}=
p_j^{K,*}-p_{j+1}^{K,*}
-\tau\partial_x H(\gamma_j^{K,*},s_j,p_j^{K,*}).
\end{equation}
Since $\eta=0$ is a right-sided minimum of $\eta\mapsto\AcK(\gamma^\eta,p^{\gamma^\eta})$, its right derivative is  nonnegative, while \eqref{eq:first-var} and (H3) give
\[
\left.\frac{\dd}{\dd\eta}
\AcK(\gamma^\eta,p^{\gamma^\eta})\right|_{\eta=0+}
=-C<0,
\]
a contradiction.

If instead $p_j^{K,*}=C$ and $p_{j+1}^{K,*}=0$, every sufficiently small $\eta<0$ makes $\gamma^\eta$ admissible. 
Since $\eta=0$ is a left-sided minimum of $\eta\mapsto\AcK(\gamma^\eta,p^{\gamma^\eta})$, its left derivative 
is nonpositive, whereas \eqref{eq:first-var} and $\partial_x H(x,s,C)=0$ give
\[
\left.\frac{\dd}{\dd\eta}
\AcK(\gamma^\eta,p^{\gamma^\eta})\right|_{\eta=0-}
=C>0,
\]
a contradiction.
\end{proof}

\begin{corollary}[Existence of an interior momentum component]\label{cor:one-interior}
If $y\in\operatorname{int}R_{t,x}^K$, then some $p_j^{K,*}$ belongs to $(0,C)$.
\end{corollary}

\begin{proof}
If every $p_k^{K,*}$ belongs to $\{0,C\}$, Lemma~\ref{lem:no-switch} forces all of them to be $0$ or all of them to be $C$.  By \eqref{eq:disc-velocity}, the corresponding curve is exactly the boundary sequence \eqref{eq:disc-ext-0} or \eqref{eq:disc-ext-C}.  Its initial point is therefore an endpoint of $R_{t,x}^K$. This contradicts the choice of $y$. 
\end{proof}

\subsection{Interiority of momentum and the discrete Hamilton's equations}

We next prove that, for each interior reachable point  $y\in {\rm int} R^K_{t,x}$ and a discrete optimizing trajectory $(\gamma^{K,*},p^{K,*})$ for the left-hand side of \eqref{discrete-inf-sup}, the momentum part $p^{K,*}$ belongs to $(0,C)^K$. 
This allows us to take two-sided variations of $\gamma^{K,*}$ to derive the discrete Hamilton's equations.

\begin{proposition}\label{prop:disc-interiority}
Let $y\in\operatorname{int}R_{t,x}^K$ and let $K$ be sufficiently large.  Then every discrete optimizing trajectory $(\gamma^{K,*},p^{K,*})$ for the left-hand side of \eqref{discrete-inf-sup} satisfies
\[
0<p_k^{K,*}<C,
\qquad k=1,\ldots,K.
\]
Consequently,
\begin{align}
\dot{\gamma}_k^{K,*}
&=
\partial_p H(\gamma_k^{K,*},s_k,p_k^{K,*})
+\sqrt\tau(p_k^{K,*}-c),
&&k=1,\ldots,K,
\label{eq:disc-Ham-x}\\
\frac{p_{k+1}^{K,*}-p_k^{K,*}}{\tau}
&=
-\partial_x H(\gamma_k^{K,*},s_k,p_k^{K,*}),
&&k=1,\ldots,K-1.
\label{eq:disc-Ham-p}
\end{align}
\end{proposition}

\begin{proof}
Corollary~\ref{cor:one-interior} gives at least one interior momentum component
$p_{j_0}^{K,*}\in(0,C)$ for some $j_0\in\{1,\ldots,K\}$.
We first prove one step of the forward propagation and then one step of the
backward propagation.

 Assume that
$p_j^{K,*}\in(0,C)$ for some $j\le K-1$.  We show that
$p_{j+1}^{K,*}\in(0,C)$.
We use the $j$-node variation $\gamma^\eta$ of $\gamma^{K,*}$: 
\[
\gamma_j^\eta:=\gamma_j^{K,*}+\eta,
\qquad
\gamma_i^\eta:=\gamma_i^{K,*}\quad(i\ne j).
\]
Then
\[
\dot{\gamma}_j^\eta
=
\dot{\gamma}_j^{K,*}+\frac{\eta}{\tau},
\qquad
\dot{\gamma}_{j+1}^\eta
=
\dot{\gamma}_{j+1}^{K,*}-\frac{\eta}{\tau}.
\]
Since $p_j^{K,*}\in(0,C)$, \eqref{eq:disc-velocity} and the strict
monotonicity in \eqref{eq:Htau-strict} show that the $j$-th velocity $\dot\gamma_j^{K,*}$ lies
strictly inside its regularized admissible interval.  Hence both signs of
sufficiently small $\eta$ preserve the $j$-th velocity constraint of $\gamma^\eta$.

Suppose first that $p_{j+1}^{K,*}=C$.  Then the $(j+1)$-st velocity lies
at its upper regularized boundary:
\[
\dot{\gamma}_{j+1}^{K,*}
=
\partial_pH(\gamma_{j+1}^{K,*},s_{j+1},C)
+\sqrt\tau(C-c).
\]
Since $\gamma_{j+1}^\eta=\gamma_{j+1}^{K,*}$ and
$\dot{\gamma}_{j+1}^\eta
=\dot{\gamma}_{j+1}^{K,*}-\eta/\tau$, we have
\[
\frac{\dd}{\dd\eta}
\left[
\partial_pH(\gamma_{j+1}^\eta,s_{j+1},C)
+\sqrt\tau(C-c)
-\dot{\gamma}_{j+1}^\eta
\right]_{\eta=0}
=
\frac1\tau>0.
\]
Hence, as discussed in the proof of Lemma~\ref{lem:no-switch}, we obtain for all sufficiently small $\eta>0$,
\[
\partial_pH(\gamma_{j+1}^\eta,s_{j+1},0)-\sqrt\tau c
<
\dot{\gamma}_{j+1}^\eta
<
\partial_pH(\gamma_{j+1}^\eta,s_{j+1},C)
+\sqrt\tau(C-c).
\]
All the other velocity constraints are unchanged.  Thus $\gamma^\eta$ is
admissible for all sufficiently small $\eta>0$, with
$\gamma_0^\eta=y$ and $\gamma_K^\eta=x$. 
We can therefore apply Lemma~\ref{lem:envelope}.  Since $\eta=0$ is a
right-sided minimum of $\eta\longmapsto\AcK(\gamma^\eta,p^{\gamma^\eta})$, 
its right derivative is nonnegative.  By \eqref{eq:first-var}, 
\[
0\le \left.\frac{\dd}{\dd\eta}\AcK(\gamma^\eta,p^{\gamma^\eta})\right|_{\eta=0+}=
p_j^{K,*}-C
-\tau\partial_xH(\gamma_j^{K,*},s_j,p_j^{K,*}).
\]
Since $\partial_xH(x,s,C)=0$ by (H3), the mean value theorem gives some
$\theta\in(p_j^{K,*},C)$ such that
\[
\partial_xH(\gamma_j^{K,*},s_j,p_j^{K,*})
=\partial_xH(\gamma_j^{K,*},s_j,p_j^{K,*})-\partial_xH(\gamma_j^{K,*},s_j,C)=
\partial_{xp}^2H(\gamma_j^{K,*},s_j,\theta)
(p_j^{K,*}-C).
\]
Consequently, in view of  \eqref{eq:tau-Hxp}, we find 
\[
0\le
(p_j^{K,*}-C)
\left(
1-\tau\partial_{xp}^2H(\gamma_j^{K,*},s_j,\theta)
\right)<0,
\]
which is a contradiction.

Suppose next that $p_{j+1}^{K,*}=0$.  Then the $(j+1)$-st velocity lies
at its lower regularized boundary:
\[
\dot{\gamma}_{j+1}^{K,*}
=
\partial_pH(\gamma_{j+1}^{K,*},s_{j+1},0)-\sqrt\tau c.
\]
For the same $j$-node variation,
\[
\frac{\dd}{\dd\eta}
\left[
\dot{\gamma}_{j+1}^\eta
-\partial_pH(\gamma_{j+1}^\eta,s_{j+1},0)
+\sqrt\tau c
\right]_{\eta=0}
=
-\frac1\tau<0.
\]
Hence, for all sufficiently small $\eta<0$, 
\[
\partial_pH(\gamma_{j+1}^\eta,s_{j+1},0)-\sqrt\tau c
<
\dot{\gamma}_{j+1}^\eta
<
\partial_pH(\gamma_{j+1}^\eta,s_{j+1},C)
+\sqrt\tau(C-c).
\]
As before, the $j$-th velocity constraint remains away from the boundary  for sufficiently
small $\eta$, and all other velocity constraints are unchanged.
Therefore $\gamma^\eta$ is admissible for all sufficiently small $\eta<0$. 
Since $\eta=0$ is a left-sided minimum of
$\eta\mapsto\AcK(\gamma^\eta,p^{\gamma^\eta})$, its left derivative is
nonpositive.  Using \eqref{eq:first-var} and
$\partial_xH(x,s,0)=0$, the mean value theorem gives some
$\theta\in(0,p_j^{K,*})$ such that
\[
0\ge\left.\frac{\dd}{\dd\eta}\AcK(\gamma^\eta,p^{\gamma^\eta})\right|_{\eta=0-}=
p_j^{K,*}
\left(
1-\tau\partial_{xp}^2H(\gamma_j^{K,*},s_j,\theta)
\right)>0,
\]
which is a contradiction.
We have therefore proved
\begin{equation}\label{eq:forward-propagation}
p_j^{K,*}\in(0,C)
\quad\Longrightarrow\quad
p_{j+1}^{K,*}\in(0,C),
\qquad j\le K-1.
\end{equation}

We next prove the backward propagation.  Assume that
$p_j^{K,*}\in(0,C)$ for some $2\le j\le K$.  We show that
$p_{j-1}^{K,*}\in(0,C)$.
Take the $(j-1)$-node variation $\gamma^\eta$ of $\gamma^{K,*}$:
\[
\gamma_{j-1}^\eta:=\gamma_{j-1}^{K,*}+\eta,
\qquad
\gamma_i^\eta:=\gamma_i^{K,*}\quad(i\ne j-1).
\]
Then
\[
\dot{\gamma}_{j-1}^\eta
=
\dot{\gamma}_{j-1}^{K,*}+\frac{\eta}{\tau},
\qquad
\dot{\gamma}_{j}^\eta
=
\dot{\gamma}_{j}^{K,*}-\frac{\eta}{\tau}.
\]
Since $p_j^{K,*}\in(0,C)$, the $j$-th velocity lies strictly inside its
regularized admissible interval.  Hence both signs of sufficiently small
$\eta$ preserve the $j$-th velocity constraint.
Note that, by Lemma~\ref{lem:envelope}, on either side for which the variation is
admissible,
\begin{equation}\label{eq:first-var-left}
\left.\frac{\dd}{\dd\eta}
\AcK(\gamma^\eta,p^{\gamma^\eta})\right|_{\eta=0\pm}
=
p_{j-1}^{K,*}-p_j^{K,*}
-\tau\partial_xH(
\gamma_{j-1}^{K,*},s_{j-1},p_{j-1}^{K,*}).
\end{equation}

Suppose first that $p_{j-1}^{K,*}=C$.  Then the $(j-1)$-st velocity lies
at its upper regularized boundary:
\[
\dot{\gamma}_{j-1}^{K,*}
=
\partial_pH(\gamma_{j-1}^{K,*},s_{j-1},C)
+\sqrt\tau(C-c).
\]
For $\eta<0$, \eqref{eq:tau-Hxp} gives
\[
\frac{\dd}{\dd\eta}
\left[
\partial_pH(\gamma_{j-1}^\eta,s_{j-1},C)
+\sqrt\tau(C-c)
-\dot{\gamma}_{j-1}^\eta
\right]_{\eta=0}
=
\partial_{xp}^2H(\gamma_{j-1}^{K,*},s_{j-1},C)
-\frac1\tau<0.
\]
Hence, for all sufficiently small $\eta<0$, 
\[
\partial_pH(\gamma_{j-1}^\eta,s_{j-1},0)-\sqrt\tau c
<
\dot{\gamma}_{j-1}^\eta
<
\partial_pH(\gamma_{j-1}^\eta,s_{j-1},C)
+\sqrt\tau(C-c).
\]
Since the $j$-th velocity constraint remains away from the boundary and all other velocity
constraints are unchanged, $\gamma^\eta$ is admissible for all sufficiently
small $\eta<0$.

Since $\eta=0$ is a left-sided minimum of
$\eta\mapsto\AcK(\gamma^\eta,p^{\gamma^\eta})$, its left derivative is
nonpositive.  However, by \eqref{eq:first-var-left} with $\partial_xH(
\gamma_{j-1}^{K,*},s_{j-1},p_{j-1}^{K,*})=\partial_xH(
\gamma_{j-1}^{K,*},s_{j-1},C)=0$ from  (H3),
\[
0\ge \left.\frac{\dd}{\dd\eta}
\AcK(\gamma^\eta,p^{\gamma^\eta})\right|_{\eta=0-}
=
C-p_j^{K,*}>0,
\]
which is a contradiction.

Suppose next that $p_{j-1}^{K,*}=0$.  Then the $(j-1)$-st velocity lies
at its lower regularized boundary:
\[
\dot{\gamma}_{j-1}^{K,*}
=
\partial_pH(\gamma_{j-1}^{K,*},s_{j-1},0)-\sqrt\tau c.
\]
For $\eta>0$, \eqref{eq:tau-Hxp} gives
\[
\frac{\dd}{\dd\eta}
\left[
\dot{\gamma}_{j-1}^\eta
-\partial_pH(\gamma_{j-1}^\eta,s_{j-1},0)
+\sqrt\tau c
\right]_{\eta=0}
=
\frac1\tau
-\partial_{xp}^2H(\gamma_{j-1}^{K,*},s_{j-1},0)>0.
\]
Hence, for all sufficiently small $\eta>0$, 
\[
\partial_pH(\gamma_{j-1}^\eta,s_{j-1},0)-\sqrt\tau c
<
\dot{\gamma}_{j-1}^\eta
<
\partial_pH(\gamma_{j-1}^\eta,s_{j-1},C)
+\sqrt\tau(C-c).
\]
Since the $j$-th velocity constraint remains away from the boundary  and all other velocity
constraints are unchanged, $\gamma^\eta$ is admissible for all sufficiently
small $\eta>0$.

Since $\eta=0$ is a right-sided minimum of
$\eta\mapsto\AcK(\gamma^\eta,p^{\gamma^\eta})$, its right derivative is
nonnegative.  On the other hand, by \eqref{eq:first-var-left} and (H3),
\[
0\le \left.\frac{\dd}{\dd\eta}
\AcK(\gamma^\eta,p^{\gamma^\eta})\right|_{\eta=0+}
=
-p_j^{K,*}<0,
\]
a contradiction.
We have therefore proved
\begin{equation}\label{eq:backward-propagation}
p_j^{K,*}\in(0,C)
\quad\Longrightarrow\quad
p_{j-1}^{K,*}\in(0,C),
\qquad 2\le j\le K.
\end{equation}

Starting from $p_{j_0}^{K,*}\in(0,C)$, we repeatedly apply
\eqref{eq:forward-propagation} if $j_0=1$, both
\eqref{eq:forward-propagation} and \eqref{eq:backward-propagation} if
$1<j_0<K$, and \eqref{eq:backward-propagation} if $j_0=K$.
Thus, in every case,
\[
0<p_k^{K,*}<C,
\qquad k=1,\ldots,K.
\]

It remains to derive the discrete Hamilton equations.  For each
$j=1,\ldots,K-1$, both $p_j^{K,*}$ and $p_{j+1}^{K,*}$ belong to
$(0,C)$.  Hence the $j$-th and $(j+1)$-st velocity constraints are both
strict, and both signs of a sufficiently small $j$-node variation are
admissible.  Therefore $\eta=0$ is a two-sided minimum of
 $\eta\longmapsto\AcK(\gamma^\eta,p^{\gamma^\eta})$. 
By Lemma~\ref{lem:envelope}, its derivative vanishes at $\eta=0$, and
\eqref{eq:first-var} gives
\[
0=p_j^{K,*}-p_{j+1}^{K,*}
-\tau\partial_xH(\gamma_j^{K,*},s_j,p_j^{K,*}).
\]
This is \eqref{eq:disc-Ham-p}.  Equation \eqref{eq:disc-Ham-x} follows
from \eqref{eq:disc-velocity}.
\end{proof}

\section{Passage to the limit $K\to\infty$ and optimality of the limit}\label{sec:optimality}

Fix $c\in(0,C)$, $t>0$, $x\in\R$, and an arbitrary reachable point
$y\in R_{t,x}$.  
By \eqref{eq:reachable-strict-inclusion}
in Lemma~\ref{lem:disc-reachable},
for all sufficiently large $K$, 
\[
y\in\operatorname{int}R_{t,x}^K.
\]
For each such $K$, let $(\gamma^{K,*},p^{K,*})$ be a discrete optimizing
trajectory for \eqref{discrete-inf-sup}.  Proposition~\ref{prop:disc-interiority} gives
\[
0<p_k^{K,*}<C,
\qquad k=1,\ldots,K,
\]
together with the discrete Hamilton equations
\eqref{eq:disc-Ham-x}--\eqref{eq:disc-Ham-p}. 
We first take the limit of these discrete Hamiltonian trajectories as  $K\to\infty$. 
We then prove that the limiting trajectory is optimal for the original
continuous two-fixed-endpoint problem.

\subsection{Linear interpolation of the discrete optimizers and its limit}
\label{subsec:Hamiltonian-limit}

Let $\gamma^K(\cdot)$ be the piecewise linear interpolation of the nodes
$\gamma_k^{K,*}$.  Let $p^K(\cdot)$ be the piecewise linear interpolation
of $p_k^{K,*}$ on $[\tau,t]$, extended constantly on $[0,\tau]$.
From \eqref{eq:disc-Ham-x} and \eqref{eq:disc-Ham-p}, 
\[
\|\dot\gamma^K(\cdot)\|_{L^\infty(0,t)}
\le
\|\partial_pH\|_{L^\infty}+\sqrt\tau C,\qquad \|\dot p^K(\cdot)\|_{L^\infty(0,t)}
\le
\|\partial_xH\|_{L^\infty}.
\]
Together with $0\le p^K\le C$ and the fixed endpoints
$\gamma^K(0)=y$, $\gamma^K(t)=x$, these estimates give uniform
boundedness and equicontinuity of $\{(\gamma^K(\cdot),p^K(\cdot))\}_{K\in\mathbb N}$.  Hence, after passing to a subsequence,
not relabeled, we find a Lipschitz trajectory $(\gamma^*(\cdot),p^*(\cdot))$ such that 
\begin{equation}\label{eq:uniform-conv}
\gamma^K(\cdot)\to\gamma^*(\cdot),
\qquad
p^K(\cdot)\to p^*(\cdot)
\qquad\text{uniformly on }[0,t]
\quad\text{as }K\to\infty.
\end{equation}
We pass to the limit in the discrete Hamilton equations in their
integral form. For each $s\in(0,t]$, take $m_s=m_s(K)\in\mathbb N$ such that $s\in 
 (m_s\tau-\tau,m_s\tau]$. 
From \eqref{eq:disc-Ham-x}, we have
\begin{align*}
 \gamma^K(s)
 &= \gamma_{m_s}^{K,*}+O(\tau)
 = y+\tau\sum_{k=1}^{m_s}
 \left\{
 \partial_pH(\gamma_k^{K,*},s_k,p_k^{K,*})
 +\sqrt\tau(p_k^{K,*}-c)
 \right\}
 +O(\tau)\\
 &=
 y+\int_0^{s_{m_s}}
 \left\{
 \partial_pH(\gamma^K(r),r,p^K(r))
 +\sqrt\tau(p^K(r)-c)
 \right\}\dd r
 +O(\tau)\\
 &=
 y+\int_0^s
 \left\{
 \partial_pH(\gamma^K(r),r,p^K(r))
 +\sqrt\tau(p^K(r)-c)
 \right\}\dd r
 +O(\tau).
\end{align*}
By \eqref{eq:uniform-conv},
\begin{equation*}
\gamma^*(s)=y+\int_0^s
 \left\{
 \partial_pH(\gamma^*(r),r,p^*(r))
 \right\}\dd r.
\end{equation*}
Similarly, we obtain 
\begin{equation*}
p^*(s)=p^\ast(0)-\int_0^s
 \left\{
 \partial_xH(\gamma^*(r),r,p^*(r))
 \right\}\dd r,\qquad p^\ast(0)\in[0,C].
\end{equation*}
This implies that $(\gamma^*(\cdot),p^\ast(\cdot))$ belongs to $C^2([0,t];\R\times[0,C])$ and satisfies  
\begin{equation*}%\label{eq:limit-Ham}
 \dot\gamma^*(s)
 =\partial_pH(\gamma^*(s),s,p^*(s)),
 \qquad
 \dot p^*(s)
 = -\partial_xH(\gamma^*(s),s,p^*(s)), 
\end{equation*}
as well as 
$\gamma^*(0)=y,
\quad
\gamma^*(t)=x, \quad 
\partial_pH(\gamma^*(s),s,0)
\le\dot\gamma^*(s)
\le\partial_pH(\gamma^*(s),s,C)$. 
So $\gamma^*(\cdot)\in\Gamma_t$. 
Since $\partial_pH(\gamma^*(s),s,p^*(s))
=\dot\gamma^*(s)$,   $p^*(s)$ is
a pointwise maximizing momentum for each $s\in[0,t]$\footnote{For every $p\in[0,C]$, convexity gives
$H(\gamma^*(s),s,p)-H(\gamma^*(s),s,p^*(s))
\ge
\partial_pH(\gamma^*(s),s,p^*(s))(p-p^*(s))
=
\dot\gamma^*(s)(p-p^*(s))$.
Thus
$\dot\gamma^*(s)(p-c)-H(\gamma^*(s),s,p)
\le
\dot\gamma^*(s)(p^*(s)-c)-H(\gamma^*(s),s,p^*(s))$.}.  Hence $p^*(\cdot)$ is a maximizing
momentum curve associated with $\gamma^*(\cdot)$, and
\begin{equation}\label{eq:AequalsF-limit}
\sup_{p(\cdot)\in P_t}
\Ac(\gamma^*(\cdot),p(\cdot))
=
\Ac(\gamma^*(\cdot),p^*(\cdot)).
\end{equation}
Moreover, by \eqref{eq:disc-Ham-x},
\begin{align*}
 \AcK(\gamma^{K,*},p^{K,*})
 &=
 \tau\sum_{k=1}^K
 \Bigl\{
 \partial_pH(\gamma_k^{K,*},s_k,p_k^{K,*})
 (p_k^{K,*}-c)
 -H(\gamma_k^{K,*},s_k,p_k^{K,*})
 +\frac{\sqrt\tau}{2}(p_k^{K,*}-c)^2
 \Bigr\}.
\end{align*}
Hence, by \eqref{eq:uniform-conv} and the same Riemann-sum argument,
\begin{align}\label{eq:min-action-conv}
 \inf_{\substack{\gamma:\,\text{admissible}\\
                  \gamma_0=y,\ \gamma_K=x}}
\left\{ \max_{p\in[0,C]^K}\AcK(\gamma,p)\right\}
 = \AcK(\gamma^{K,*},p^{K,*})
 \longrightarrow
 \Ac(\gamma^*(\cdot),p^*(\cdot))
\quad \mbox{as $K\to\infty$}.
\end{align}

\subsection{Optimality of the limiting Hamiltonian trajectory}
\label{subsec:optimality-limit}

It remains to prove that the Hamiltonian trajectory obtained in
Subsection~\ref{subsec:Hamiltonian-limit} is optimal for the original
two-fixed-endpoint  $\inf$-$\sup$ problem.  For this purpose, we discretize an
arbitrary admissible curve.

\begin{lemma}
\label{lem:recovery}
Let $\gamma(\cdot)\in\Gamma_t$ satisfy
$\gamma(0)=y$ and $\gamma(t)=x$.  For $s_k=k\tau$, define
$\gamma^K\in\R^{K+1}$ by
\begin{equation*}
 \gamma_k^K:=\gamma(s_k),
 \qquad
 \dot\gamma_k^K
 :=
 \frac{\gamma(s_k)-\gamma(s_{k-1})}{\tau}.
\end{equation*}
Then, for all sufficiently large $K$, the discrete curve $\gamma^K$
is admissible for the regularized discrete problem
\eqref{eq:disc-adm}.  Moreover,
\begin{equation}\label{eq:recovery-conv}
 \max_{p\in[0,C]^K}\AcK(\gamma^K,p)
 \longrightarrow
 \sup_{p(\cdot)\in P_t}
 \Ac(\gamma(\cdot),p(\cdot))
 \qquad\text{as }K\to\infty.
\end{equation}
\end{lemma}

\begin{proof}
Since $\gamma(\cdot)\in\Gamma_t$, there is a constant $M_0>0$ such that
$|\dot\gamma(s)|\le M_0$ for almost every $s$.
%Let $L_0$ be a common Lipschitz constant in $(x,s)$ for $\partial_pH(x,s,0)$ and $\partial_pH(x,s,C)$.
For any $s\in[s_{k-1},s_k]$,
\begin{equation*}
 |\gamma(s)-\gamma(s_k)|\le M_0\tau,
 \qquad
 |s-s_k|\le\tau.
\end{equation*}
Hence, for almost every $s\in[s_{k-1},s_k]$, we have 
\begin{align*}
 \dot\gamma(s)
 &\ge
 \partial_pH(\gamma(s),s,0)
 \ge
 \partial_pH(\gamma(s_k),s_k,0)-D\tau,\\
 \dot\gamma(s)
 &\le
 \partial_pH(\gamma(s),s,C)
 \le
 \partial_pH(\gamma(s_k),s_k,C)+D\tau
\end{align*}
for some constant $D>0$ independent of $k$ and $K$.
Averaging over $[s_{k-1},s_k]$ gives
\begin{equation}\label{eq:average-constraint}
 \partial_pH(\gamma_k^K,s_k,0)-D\tau
 \le
 \dot\gamma_k^K
 \le
 \partial_pH(\gamma_k^K,s_k,C)+D\tau.
\end{equation}
Since $c>0$ and $C-c>0$, for all sufficiently small $\tau$,
\begin{equation*}
 D\tau
 \le
 \sqrt\tau\min\{c,C-c\}.
\end{equation*}
Thus \eqref{eq:average-constraint} implies the discrete admissibility
condition \eqref{eq:disc-adm}\footnote{This is one of the reasons why we regularize $H$ with $\frac{\sqrt{\tau}}{2}(p-c)^2$ instead of $\frac{\sqrt{\tau}}{2}p^2$}. 

We next compare the maximized actions.  Observe that 
%the maximization with respect to $p=(p_1,\ldots,p_K)\in[0,C]^K$ is performed independently on each cell.  Hence
\begin{align*}
 \max_{p\in[0,C]^K}\AcK(\gamma^K,p)
 &=
 \sum_{k=1}^K
 \tau
 \max_{p\in[0,C]}
 \left\{
 \dot\gamma_k^K(p-c)
 -H(\gamma(s_k),s_k,p)
 -\frac{\sqrt\tau}{2}(p-c)^2
 \right\}\\
 &=
 \sum_{k=1}^K
 \int_{s_{k-1}}^{s_k}
 \max_{p\in[0,C]}
 \left\{
 \dot\gamma_k^K(p-c)
 -H(\gamma(s_k),s_k,p)
 -\frac{\sqrt\tau}{2}(p-c)^2
 \right\}\dd s,
\end{align*}
while \eqref{eq:sup-attained} in Proposition~\ref{prop:measurable-max} gives
\begin{align*}
 \sup_{p(\cdot)\in P_t}
 \Ac(\gamma(\cdot),p(\cdot))
 &=
 \int_0^t
 \max_{0\le p\le C}
 \left\{
 \dot\gamma(s)(p-c)
 -
 H(\gamma(s),s,p)
 \right\}\dd s\\
 &=
 \sum_{k=1}^K
 \int_{s_{k-1}}^{s_k}
 \max_{p\in[0,C]}
 \left\{
 \dot\gamma(s)(p-c)
 -
 H(\gamma(s),s,p)
 \right\}\dd s.
\end{align*}
Therefore,
\begin{align*}
&\left|
 \max_{p\in[0,C]^K}\AcK(\gamma^K,p)
 -
 \sup_{p(\cdot)\in P_t}\Ac(\gamma(\cdot),p(\cdot))
\right|\\
&\le
 \sum_{k=1}^K
 \int_{s_{k-1}}^{s_k}
 \left|
 \max_{p\in[0,C]}
 \left\{
 \dot\gamma_k^K(p-c)
 -
 H(\gamma(s_k),s_k,p)
 -
 \frac{\sqrt\tau}{2}(p-c)^2
 \right\}
 \right.\\
&\hspace{43mm}\left.
 -
 \max_{p\in[0,C]}
 \left\{
 \dot\gamma(s)(p-c)
 -
 H(\gamma(s),s,p)
 \right\}
 \right|
 \dd s\\
&\le
 \sum_{k=1}^K
 \int_{s_{k-1}}^{s_k}
 \max_{p\in[0,C]}
 \Biggl|
 (\dot\gamma_k^K-\dot\gamma(s))(p-c)
 -H(\gamma(s_k),s_k,p)
 +H(\gamma(s),s,p)
 - \frac{\sqrt\tau}{2}(p-c)^2
 \Biggr|\dd s.
\end{align*}
Since $c\in(0,C)$ and $p\in[0,C]$, we have $|p-c|\le C$.
Also, if $K_H$ is a Lipschitz constant of $H$ in $(x,s)$ on the
compact set under consideration, then
\begin{equation*}
 |H(\gamma(s_k),s_k,p)-H(\gamma(s),s,p)|
 \le
 K_H
 \left\{
 |\gamma(s_k)-\gamma(s)|
 +
 |s_k-s|
 \right\}.
\end{equation*}
Consequently,
\begin{align*}
&\left|
 \max_{p\in[0,C]^K}\AcK(\gamma^K,p)
 - \sup_{p(\cdot)\in P_t}\Ac(\gamma(\cdot),p(\cdot))
\right|\\
&\le
 C\sum_{k=1}^K
 \int_{s_{k-1}}^{s_k}
 |\dot\gamma_k^K-\dot\gamma(s)|\dd s
 +K_H\sum_{k=1}^K
 \int_{s_{k-1}}^{s_k}
 \left\{|\gamma(s_k)-\gamma(s)|
 + |s_k-s|
 \right\}\dd s
 + \frac{t\sqrt\tau}{2}C^2.
\end{align*}
Extending $\dot\gamma$ by zero outside $(0,t)$, we have 
\begin{align*}
&\sum_{k=1}^K
 \int_{s_{k-1}}^{s_k}
 |\dot\gamma_k^K-\dot\gamma(s)|\dd s
\le
 \frac{1}{\tau}
 \sum_{k=1}^K
 \int_{s_{k-1}}^{s_k}
 \int_{s_{k-1}}^{s_k}
 |\dot\gamma(r)-\dot\gamma(s)|\dd r\dd s \\
 &\le
 \frac{1}{\tau}
 \int_{-\tau}^{\tau}
 \int_{\mathbb R}
 |\dot\gamma(s+h)-\dot\gamma(s)|\dd s\dd h
\le
 2\sup_{|h|\le\tau}
 \int_{\mathbb R}
 |\dot\gamma(s+h)-\dot\gamma(s)|\dd s
 \to0\quad\text{as }K\to\infty. 
\end{align*}
Since $|\gamma(s_k)-\gamma(s)|\le M_0|s_k-s|$, we have 
\begin{align*}
&\sum_{k=1}^K
 \int_{s_{k-1}}^{s_k}
 \left\{
 |\gamma(s_k)-\gamma(s)|
 + |s_k-s|
 \right\}\dd s
 \le(M_0+1)
 \sum_{k=1}^K
 \int_{s_{k-1}}^{s_k}
 |s_k-s|\dd s
= \frac{t(M_0+1)}{2}\tau.
\end{align*}
Now \eqref{eq:recovery-conv} is clear. 
\end{proof}

\begin{theorem}[Existence of an optimal Hamiltonian trajectory]
\label{thm:fixed-optimal}
Assume (H1)--(H3) and fix $c\in(0,C)$, $t>0$, $x\in\R$, and
$y\in R_{t,x}$.  Then there exist
$\gamma^*\in C^2([0,t];\R)$ and
$p^*\in C^2([0,t];[0,C])$ such that $\gamma^*(0)=y$, $\gamma^*(t)=x$, 
\begin{align}
 \Ac(\gamma^*(\cdot),p^*(\cdot))
 &=
 \sup_{p(\cdot)\in P_t}
 \Ac(\gamma^*(\cdot),p(\cdot))
 =
 \inf_{\substack{\gamma(\cdot)\in\Gamma_t\\
                  \gamma(0)=y,\ \gamma(t)=x}}
 \sup_{p(\cdot)\in P_t}
 \Ac(\gamma(\cdot),p(\cdot)),
 \label{eq:thm-optimal}
\end{align}
and $\dot\gamma^*(s)
 = \partial_pH(\gamma^*(s),s,p^*(s))$, 
$\dot p^*(s)
 = -\partial_xH(\gamma^*(s),s,p^*(s))$.
%Thus $(\gamma^*(\cdot),p^*(\cdot))$ is a $C^2$-smooth optimal trajectory satisfying Hamilton's equations.
\end{theorem}

\begin{proof}
Let $(\gamma^*(\cdot),p^*(\cdot))$ be the Hamiltonian trajectory
obtained in Subsection~\ref{subsec:Hamiltonian-limit}.
By \eqref{eq:AequalsF-limit}, $p^*(\cdot)$ is a maximizing momentum
curve associated with $\gamma^*(\cdot)$.

Let $\gamma(\cdot)\in\Gamma_t$ be any admissible curve satisfying
$\gamma(0)=y$ and $\gamma(t)=x$.
By Lemma~\ref{lem:recovery}, its discretization $\gamma^K$ is
discrete admissible for all sufficiently large $K$, and
\begin{equation*}
 \max_{p\in[0,C]^K}\AcK(\gamma^K,p)
 \longrightarrow
 \sup_{p(\cdot)\in P_t}
 \Ac(\gamma(\cdot),p(\cdot)) \qquad\text{as }K\to\infty.
\end{equation*}
Since $(\gamma^{K,*},p^{K,\ast})$ is a discrete optimizer, we have 
\begin{equation*}
 \AcK(\gamma^{K,*},p^{K,*})
 \le
 \max_{p\in[0,C]^K}\AcK(\gamma^K,p).
\end{equation*}
Passing to the limit and using \eqref{eq:min-action-conv} and \eqref{eq:recovery-conv}, we obtain
\begin{equation*}
 \Ac(\gamma^*(\cdot),p^*(\cdot))=\sup_{p(\cdot)\in P_t}
 \Ac(\gamma^*(\cdot),p(\cdot))
 \le
 \sup_{p(\cdot)\in P_t}
 \Ac(\gamma(\cdot),p(\cdot)).
\end{equation*}
Since $\gamma(\cdot)$ is arbitrary, we conclude 
\eqref{eq:thm-optimal}. \end{proof}

\begin{remark}
\label{rem:selection-not-all}
Even for a $C^2$-smooth minimizing configuration-space curve,
Proposition~\ref{prop:measurable-max} provides an associated measurable
maximizing momentum curve, which may differ from the $C^2$-smooth
Hamiltonian momentum selected in Theorem~\ref{thm:fixed-optimal},
particularly in the degenerate case. 
Thus the variational definition of optimality does not, by itself, guarantee that an optimal trajectory is a classical Hamiltonian trajectory. 
It is also not clear whether there exist other minimizing
configuration-space curves having less regularity.  This is related to
the weak compactness difficulty described in the Introduction.
\end{remark}
\begin{corollary}[Interior and boundary reachable points]
\label{cor:reachable-optimal}
Under the assumptions of Theorem~\ref{thm:fixed-optimal} with  $R_{t,x}=[\overline{\gamma}(0),\underline{\gamma}(0)]$, the following
statements hold.
\begin{enumerate}[label=\textup{(\roman*)}]
\item
If $y\in\operatorname{int}R_{t,x}$, the optimal Hamiltonian trajectory
obtained above satisfies
\begin{equation*}
 0<p^*(s)<C
 \qquad\text{for every }s\in[0,t].
\end{equation*}

\item
If $y=\overline{\gamma}(0)$, then for every optimal trajectory $(\gamma^\ast(\cdot),p^\ast(\cdot))$ connecting $x$ and $y$ we have  $\gamma^\ast(\cdot)=\overline{\gamma}(\cdot)$. Furthermore, 
$(\overline{\gamma}(\cdot),C)$ is an optimal Hamiltonian trajectory connecting $x$ and $y$.

\item
If $y=\underline{\gamma}(0)$, then for every optimal trajectory $(\gamma^\ast(\cdot),p^\ast(\cdot))$ connecting $x$ and $y$ we have  $\gamma^\ast(\cdot)=\underline{\gamma}(\cdot)$. Furthermore, 
$(\underline{\gamma}(\cdot),0)$ is an optimal Hamiltonian trajectory  connecting $x$ and $y$.
\end{enumerate}
\end{corollary}

\begin{proof}
Recall that by Proposition~\ref{prop:reachable} an admissible configuration-space curve $\gamma(\cdot)$ connecting $x$ and  $\underline{\gamma}(0)$ (resp. $\overline{\gamma}(0)$) is only $\underline{\gamma}(\cdot)$ (resp. $\overline{\gamma}(\cdot)$);   
 $(\overline{\gamma}(\cdot),C)$ and $(\underline{\gamma}(\cdot),0)$ satisfy Hamilton's equations. By (H3), $\R\times\{p=0\}$ and $\R\times\{p=C\}$ are invariant under Hamilton's equations.   

(i): Suppose that $p^*(s_0)\in\{0,C\}$ for some $s_0\in[0,t]$. 
Then,  the invariance implies that $(\gamma^*(\cdot),p^*(\cdot))=(\gamma^*(\cdot), 0)$ or $(\gamma^*(\cdot), C)$ with $\gamma^*(t)=x$. Hence,  the uniqueness of \eqref{eq:extreme-C} and \eqref{eq:extreme-0} implies that $\gamma^*(\cdot)$ coincides with $\underline{\gamma}(\cdot)$ or $\overline{\gamma}(\cdot)$.   Proposition~\ref{prop:reachable}  gives $\gamma^\ast(0)=y\in\{\overline{\gamma}(0),\underline{\gamma}(0)\}$, which contradicts $y\in\operatorname{int}R_{t,x}$.
Therefore,  $0<p^*(s)<C$ for every $s\in[0,t]$. 

(ii): By uniqueness of an admissible configuration-space curve connecting $x$ and $y$, we have $\gamma^*(\cdot) = \overline{\gamma}(\cdot)$; furthermore,   together with $\dot{\overline{\gamma}}(s)= \partial_pH({\overline{\gamma}}(s),s,C)$, we find that  
\begin{align*}
& \inf_{\substack{\gamma(\cdot)\in\Gamma_t\\
                  \gamma(0)=y,\ \gamma(t)=x}}\sup_{p(\cdot)\in P_t}\Ac(\gamma(\cdot),p(\cdot))
= \sup_{p(\cdot)\in P_t}\Ac(\overline{\gamma}(\cdot),p(\cdot))    \\
&\quad =\int_0^t\max_{p\in[0,C]}\{\dot{\overline{\gamma}}(s)(p-c)-H(\overline{\gamma}(s),s,p)\}\dd s 
=       \Ac(\overline{\gamma},C)         
\end{align*}
This shows the assertion. 

The parallel argument proves (iii). 
\end{proof}

\section{The Bolza  inf-sup problem}\label{sec:bolza}

We consider the corresponding Bolza problem, that is, the
problem with one fixed endpoint and an initial cost.
Let $w\in C(\T;\R)$, identified with its $2\pi$-periodic lift to $\R$.
For fixed $c\in(0,C)$, $t>0$, and $x\in\R$, we consider
\begin{equation}\label{eq:Bolza-problem}
 \inf_{\substack{\gamma(\cdot)\in\Gamma_t\\
                  \gamma(t)=x}}
 \left\{
 \sup_{p(\cdot)\in P_t}
 \Ac(\gamma(\cdot),p(\cdot))
 +
 w(\gamma(0))
 \right\}.
\end{equation}

\begin{theorem}[Existence of an optimal Hamiltonian trajectory]
\label{thm:Bolza}
Assume (H1)--(H3), and fix $c\in(0,C)$, $t>0$, $x\in\R$, and
$w\in C(\T;\R)$.  Then there exist
$\gamma^*\in C^2([0,t];\R)$ and
$p^*\in C^2([0,t];[0,C])$ such that
$\gamma^*(t)=x$,
\begin{align*}
 \Ac(\gamma^*(\cdot),p^*(\cdot))
 +w(\gamma^*(0))\notag
&= \sup_{p(\cdot)\in P_t}
 \Ac(\gamma^*(\cdot),p(\cdot))
 + w(\gamma^*(0))\notag\\
 &= \inf_{\substack{\gamma(\cdot)\in\Gamma_t\\
                  \gamma(t)=x}}
 \left\{
 \sup_{p(\cdot)\in P_t}
 \Ac(\gamma(\cdot),p(\cdot))
 + w(\gamma(0))
 \right\}, %\label{eq:Bolza-opt}
\end{align*}
and  $\dot\gamma^*(s)
 = \partial_pH(\gamma^*(s),s,p^*(s))$, 
$\dot p^*(s)
 = -\partial_xH(\gamma^*(s),s,p^*(s))$. 
\end{theorem}

\begin{proof}
We first note that the free initial endpoint is trapped in the bounded
set $R_{t,x}$ by the admissibility condition. 
Since $w$ is continuous, its contribution is uniformly
bounded and the Bolza infimum is finite.

Let $\gamma_n(\cdot)\in\Gamma_t$, $\gamma_n(t)=x$, be a minimizing
sequence for \eqref{eq:Bolza-problem}, and set
\begin{equation*}
 y_n:=\gamma_n(0).
\end{equation*}
Then $y_n\in R_{t,x}$ for every $n$.
For each $n$, Theorem~\ref{thm:fixed-optimal} gives an optimal
Hamiltonian trajectory $(\gamma_n^*(\cdot),p_n^*(\cdot))$ such that $\gamma_n^*(0)=y_n$, $\gamma_n^*(t)=x$ and 
\begin{align*}
 \Ac(\gamma_n^*(\cdot),p_n^*(\cdot))
 &=
 \inf_{\substack{\gamma(\cdot)\in\Gamma_t\\
                  \gamma(0)=y_n,\ \gamma(t)=x}}
\left\{ \sup_{p(\cdot)\in P_t}
 \Ac(\gamma(\cdot),p(\cdot))\right\}
 \le \sup_{p(\cdot)\in P_t}
 \Ac(\gamma_n(\cdot),p(\cdot)).
\end{align*}
Since $\gamma_n^*(0)=\gamma_n(0)=y_n$, the initial cost is unchanged.
Hence
\begin{align*}
 &\Ac(\gamma_n^*(\cdot),p_n^*(\cdot))+w(y_n)
 \le \sup_{p(\cdot)\in P_t}
 \Ac(\gamma_n(\cdot),p(\cdot))
 +w(\gamma_n(0)).
\end{align*}
On the other hand,
$(\gamma_n^*(\cdot),p_n^*(\cdot))$ is itself an admissible trajectory 
for the Bolza problem.  Therefore,
\begin{align*}
 &\inf_{\substack{\gamma(\cdot)\in\Gamma_t\\
                  \gamma(t)=x}}
 \left\{
 \sup_{p(\cdot)\in P_t}
 \Ac(\gamma(\cdot),p(\cdot))
 +w(\gamma(0))
 \right\}
\le \Ac(\gamma_n^*(\cdot),p_n^*(\cdot))+w(y_n).
\end{align*}
Since $\gamma_n(\cdot)$ is a minimizing sequence, it holds  that
\begin{align*}
 \Ac(\gamma_n^*(\cdot),p_n^*(\cdot))+w(y_n)
\longrightarrow
 \inf_{\substack{\gamma(\cdot)\in\Gamma_t\\
                  \gamma(t)=x}}
 \left\{
 \sup_{p(\cdot)\in P_t}
 \Ac(\gamma(\cdot),p(\cdot))
 +w(\gamma(0))
 \right\} \qquad \mbox{ as $n\to\infty$}.
\end{align*}
%Thus the Bolza minimizing sequence may be chosen to consist of classical Hamiltonian trajectories.

Since $y_n\in R_{t,x}$ and $R_{t,x}$ is compact, and $(\gamma_n^*(\cdot),p_n^*(\cdot))$ is uniformly bounded and equicontinuous, after passing to a
subsequence, we find $y^\ast$, $\gamma^*(\cdot)$ and  $p^*(\cdot)$ such that 
\begin{align*}
& y_n\longrightarrow y^*\in R_{t,x}\qquad \mbox{as $n\to\infty$},\\
&  \gamma_n^*(\cdot)\to\gamma^*(\cdot),
 \quad
 p_n^*(\cdot)\to p^*(\cdot)
 \qquad\text{uniformly on }[0,t],\qquad \mbox{as $n\to\infty$}.
\end{align*}
It follows from the integral form of Hamilton's equations  that $(\gamma^\ast(\cdot),p^\ast(\cdot))$ is a $C^2$-smooth trajectory of Hamilton's equations with $\gamma^\ast(0)=y^\ast$, $\gamma^\ast(t)=x$; this leads to   $\gamma^*(\cdot)\in\Gamma_t$, $0\le p^\ast(s)\le C$ for all $s\in[0,t]$ and 
$$ 
\Ac(\gamma^*(\cdot),p^*(\cdot))=\sup_{p(\cdot)\in P_t}\Ac(\gamma^*(\cdot),p(\cdot)),$$
i.e.,   $p^*(\cdot)$ is a maximizing momentum curve associated with
$\gamma^*(\cdot)$. 
Since each $(\gamma_n^*(\cdot),p_n^*(\cdot))$ satisfies Hamilton's
equations, we have 
\begin{align*}
 \Ac(\gamma_n^*(\cdot),p_n^*(\cdot))
 &=
 \int_0^t
 \left\{
 \partial_pH(\gamma_n^*(s),s,p_n^*(s))
 (p_n^*(s)-c)
 - H(\gamma_n^*(s),s,p_n^*(s))
 \right\}\dd s,
\end{align*}
from which we obtain 
\begin{equation*}
 \Ac(\gamma_n^*(\cdot),p_n^*(\cdot))
 \longrightarrow
 \Ac(\gamma^*(\cdot),p^*(\cdot))\qquad\mbox{ as $n\to\infty$}. 
\end{equation*}
Since  $w(y_n)\to w(y^*)=w(\gamma^\ast(0))$ as $n\to\infty$, we conclude that 
\begin{align*}
 \Ac(\gamma^*(\cdot),p^*(\cdot))+w(\gamma^\ast(0))
 &=  \sup_{p(\cdot)\in P_t} \left\{
 \Ac(\gamma^\ast(\cdot),p(\cdot))
 +w(\gamma^\ast(0))
 \right\}\\
 &= \inf_{\substack{\gamma(\cdot)\in\Gamma_t\\
                  \gamma(t)=x}}
 \left\{
 \sup_{p(\cdot)\in P_t}
 \Ac(\gamma(\cdot),p(\cdot))
 +w(\gamma(0))
 \right\}.
\end{align*}
Therefore $(\gamma^*(\cdot),p^*(\cdot))$ is an optimal Hamiltonian
trajectory for the Bolza problem.
\end{proof}
%%%%%%%%%%%%%%%%%

\medskip\medskip\medskip

\noindent{\bf Acknowledgements.} 
Takaaki Nishida is supported by JSPS KAKENHI 26K06864 and partially by JSPS KAKENHI 24H00185. 
  Kohei Soga is supported by JSPS KAKENHI 22K03391 and partially by JSPS KAKENHI  26H02003.

%%%%%%%%%%%%%%
\medskip\medskip\medskip

\noindent{\bf Data availability.}
Data sharing is not applicable to this article as no datasets were generated or analyzed during the current study.

\medskip

\noindent{\bf Conflicts of interest statement.} The authors state that there is no conflict of interest.
%%%%%%%%%%%%%%%%%

\begin{appendices}
\section{Reduction of the two-degree-of-freedom system}\label{app:reduction}

We verify equivalence between the flow of $\mathscr H$ on $\{\mathscr H=h\}$ and that of $f$ mentioned in Section~\ref{sec:model-motivation}, where the argument below is independent of the specific form of $\mathscr H$ and $f$. 

Differentiating \eqref{eq:f-energy-identity}, rather than the explicit square-root formula, gives
\begin{align}
 \partial_r f(r,t,y)
 &=\frac{\partial_{x_1}\mathscr H\bigl(r,t,y,-f(r,t,y)\bigr)}
         {\partial_{p_2}\mathscr H\bigl(r,t,y,-f(r,t,y)\bigr)},\label{eq:f-r}\\
 \partial_t f(r,t,y)
 &=\frac{\partial_{x_2}\mathscr H\bigl(r,t,y,-f(r,t,y)\bigr)}
         {\partial_{p_2}\mathscr H\bigl(r,t,y,-f(r,t,y)\bigr)},\notag\\
 \partial_y f(r,t,y)
 &=\frac{\partial_{p_1}\mathscr H\bigl(r,t,y,-f(r,t,y)\bigr)}
         {\partial_{p_2}\mathscr H\bigl(r,t,y,-f(r,t,y)\bigr)}.\label{eq:f-y}
\end{align}

Let $(x_1(s),x_2(s),p_1(s),p_2(s))$ be a Hamiltonian trajectory of $\mathscr H$ contained in the energy surface $\{\mathscr H=h\}$.  By \eqref{eq:model-p2-derivative},
\[
 \dot x_2(s)
 =\partial_{p_2}\mathscr H\bigl(x_1(s),x_2(s),p_1(s),p_2(s)\bigr)
 \ge\omega>0.
\]
Hence $s\mapsto x_2(s)$ is strictly increasing.  Let $s=\eta(t)$ denote its inverse and set $r(t):=x_1(\eta(t))$ and $y(t):=p_1(\eta(t))$.  Since the trajectory lies on $\{\mathscr H=h\}$, $p_2(\eta(t))=-f(r(t),t,y(t))$, and
\[
 \eta'(t)
 =\frac{1}{\partial_{p_2}\mathscr H\bigl(r(t),t,y(t),-f(r(t),t,y(t))\bigr)}.
\]
Using Hamilton's equations for $\mathscr H$ together with \eqref{eq:f-y} and \eqref{eq:f-r}, we obtain
\begin{align*}
 r'(t)
 &=\dot x_1(\eta(t))\eta'(t)
 =\frac{\partial_{p_1}\mathscr H\bigl(r(t),t,y(t),-f(r(t),t,y(t))\bigr)}
         {\partial_{p_2}\mathscr H\bigl(r(t),t,y(t),-f(r(t),t,y(t))\bigr)}
 =\partial_y f(r(t),t,y(t)),\\
 y'(t)
 &=\dot p_1(\eta(t))\eta'(t)
 =-\frac{\partial_{x_1}\mathscr H\bigl(r(t),t,y(t),-f(r(t),t,y(t))\bigr)}
          {\partial_{p_2}\mathscr H\bigl(r(t),t,y(t),-f(r(t),t,y(t))\bigr)}
 =-\partial_r f(r(t),t,y(t)).
\end{align*}
Hence $(r(t),y(t))$ satisfies the reduced system \eqref{eq:reduced-Ham}.

Conversely, let $(r(t),y(t))$ solve \eqref{eq:reduced-Ham} with $0\le y(t)\le U(r(t);h)$.  Fix $s_0$ and $t_0$, set $x_2(s_0)=t_0$, and solve
\begin{equation}\label{eq:appendix-x2}
 \dot x_2(s)
 =\partial_{p_2}\mathscr H\Bigl(
 r(x_2(s)),x_2(s),y(x_2(s)),
 -f(r(x_2(s)),x_2(s),y(x_2(s)))\Bigr).
\end{equation}
The right-hand side is positive by \eqref{eq:model-p2-derivative}.  Define
\begin{equation}\label{eq:appendix-reconstruction}
 x_1(s):=r(x_2(s)),\qquad
 p_1(s):=y(x_2(s)),\qquad
 p_2(s):=-f(x_1(s),x_2(s),p_1(s)).
\end{equation}
By \eqref{eq:implicit-f}, the reconstructed curve stays on $\{\mathscr H=h\}$.  Moreover,
\begin{align*}
 \dot x_1(s)
 &=\partial_y f(x_1(s),x_2(s),p_1(s))\,\dot x_2(s)\\
 &=\frac{\partial_{p_1}\mathscr H\bigl(x_1(s),x_2(s),p_1(s),-f(x_1(s),x_2(s),p_1(s))\bigr)}
         {\partial_{p_2}\mathscr H\bigl(x_1(s),x_2(s),p_1(s),-f(x_1(s),x_2(s),p_1(s))\bigr)}\,\dot x_2(s)\\
 &=\partial_{p_1}\mathscr H\bigl(x_1(s),x_2(s),p_1(s),p_2(s)\bigr),\\
 \dot p_1(s)
 &=-\partial_r f(x_1(s),x_2(s),p_1(s))\,\dot x_2(s)\\
 &=-\frac{\partial_{x_1}\mathscr H\bigl(x_1(s),x_2(s),p_1(s),-f(x_1(s),x_2(s),p_1(s))\bigr)}
          {\partial_{p_2}\mathscr H\bigl(x_1(s),x_2(s),p_1(s),-f(x_1(s),x_2(s),p_1(s))\bigr)}\,\dot x_2(s)\\
 &=-\partial_{x_1}\mathscr H\bigl(x_1(s),x_2(s),p_1(s),p_2(s)\bigr).
\end{align*}
The equation \eqref{eq:appendix-x2} is exactly $\dot x_2(s)=\partial_{p_2}\mathscr H\bigl(x_1(s),x_2(s),p_1(s),p_2(s)\bigr)$.
Finally, differentiating the last identity in \eqref{eq:appendix-reconstruction} gives
\begin{align*}
 \dot p_2(s)
 =&-\partial_r f(x_1(s),x_2(s),p_1(s))\,\dot x_1(s)
 -\partial_t f(x_1(s),x_2(s),p_1(s))\,\dot x_2(s)\\
 &-\partial_y f(x_1(s),x_2(s),p_1(s))\,\dot p_1(s).
\end{align*}
With \eqref{eq:f-r}--\eqref{eq:f-y} and the three equations already obtained, the first and third terms cancel, while the second term gives
\[
 \dot p_2(s)
 =-\partial_{x_2}\mathscr H\bigl(x_1(s),x_2(s),p_1(s),p_2(s)\bigr).
\]
Thus the reconstructed curve satisfies all four of Hamilton's equations for $\mathscr H$.  This proves the equivalence stated in Section~\ref{sec:model-motivation}. 

\end{appendices}

\end{document}